\documentclass[11pt]{amsart}
\usepackage{amsmath}
\usepackage[dvipsnames]{xcolor}
\usepackage[normalem]{ulem}
\usepackage{mathrsfs}
\usepackage{amssymb}

\newtheorem{theorem}{Theorem}[section]
\newtheorem{theoremx}{Theorem}

\newtheorem{claim}[theorem]{Claim}

\newtheorem{lemma}[theorem]{Lemma}

\newtheorem{proposition}[theorem]{Proposition}
\newtheorem{corollary}[theorem]{Corollary}

\theoremstyle{definition}
\newtheorem{definition}[theorem]{Definition}

\newtheorem{question}[theorem]{Question}

\theoremstyle{remark}
\newtheorem{remark}[theorem]{Remark}

\newtheorem{notation}[theorem]{Notation}

\def\l{{\langle}}
\def\r{{\rangle}}

\DeclareMathOperator{\dom}{dom}
\DeclareMathOperator{\supp}{Supp}

\DeclareMathOperator{\rng}{rng}
\newcommand{\cf}{{\rm cf}}
\def\Ult{{\rm Ult}}

\title{A small ultrafilter number at every singular cardinal}
\author{Tom Benhamou}
\thanks{The research of the first author was supported by the National Science Foundation under Grant
No. DMS-2246703}
\address[Benhamou]{Department of Mathematics, Statistics, and Computer Science, University of Illinois at Chicago, Chicago, IL 60607, USA}
\email{tomb@uic.edu}

\author{Sittinon Jirattikansakul}
\thanks{The second author was partially supported by ISF grant No. 882/22.
The authors would like to thank Moti Gitik for a useful discussion.}
\address[Jirattikansakul]{School of Mathematical Sciences, Tel Aviv University, Tel Aviv-Yafo,  Tel Aviv, Israel, 6997801}
\email{jir.sittinon@gmail.com}

\date{February 2023}

\begin{document}
\maketitle
\begin{abstract}
    We obtain a small ultrafilter number at $\aleph_{\omega_1}$. Moreover, we develop a version of the overlapping strong extender forcing with collapses which can keep the top cardinal $\kappa$ inaccessible. We apply this forcing to construct  a model where $\kappa$ is the least inaccessible and $V_\kappa$ is a model of GCH at regulars, failures of SCH at singulars, and the ultrafilter numbers at all singulars are small. 
\end{abstract}
\maketitle
\section{introduction}
Some of the most basic mathematical theorems rely on the possibility to distinguish between \textit{small} and \textit{large} sets. For example, the Lebesgue criterion for Riemann integrability states that a bounded function $f:[a,b]\rightarrow \mathbb{R}$ is Riemann integrable, if and only if the set of its discontinuity points is \textit{small}, which in this case means of Lebesgue measure zero. Smallness has many other interpretations: small cardinalities in set theory, nowhere dense sets in topology, probability zero events in probability theory, or polynomial and linear functions in computability theory. An abstract approach to defining a notion of largeness for the subsets of a given set $X$ is a \textit{filter}, which is simply a set $F\subseteq P(X)$ that contains all the large subsets of $X$. Formally speaking, we require the following $3$ axioms which say that $F$ is a \textit{filter over $X$}:
\begin{enumerate}
    \item  $X\in F$, $\emptyset\notin F$. (non empty and non-degenerate)
    \item $A,B\in F\Rightarrow A\cap B\in F$. (closed under intersection)
    \item $(A\in F\wedge A\subseteq B)\Rightarrow B\in F$. (upward closed to inclusion)
\end{enumerate}

For a fixed filter $F$, we may consider small sets as the sets whose complements are in $F$. Note that for many filters, there are sets $X$ which are neither small, nor large, namely $X\notin F$ and also $X^c\notin F$ e.g.  in probability, there are sets $X$ such that $0<\mathbb{P}(X)<1$ and thus are neither small (i.e. $\mathbb{P}(X)=0$), nor large (i.e. $\mathbb{P}(X)=1$).

Ultrafilters are those filters which do determine that every set is either large or small. Namely, a filter $U$ over $X$ is an \textit{ultrafilter} if for every  $B\subseteq X$ either $B\in U$ or $X\setminus B\in U$. Most of the non-trivial examples of ultrafilters involve the Axiom of choice and are thus highly non-constructive. However, they have been proven to be useful in many areas such as Analysis, Topology, Model theory, Algebra, and Combinatorics. For example, Nonstandard analysis \cite{Henson1997} is an alternative approach to study analysis and more sophisticated mathematics. The  $\epsilon$-$\delta$ definitions in analysis can be replaced by using more concrete objects, so-called  infinitesimals, in a nonstandard universe of nonstandard reals, which contain all reals. One of the constructions of a nonstandard universe is through an ultraproduct construction, which requires a non-trivial ultrafilter over $\mathbb{N}$. One especially useful usage of ultrafilters appears in topology, more specifically, in the construction of the Stone-\v{C}ech compactification \cite{stone1937},\cite{cech1937}.  Studying the combinatorial nature of ultrafilters is important to obtain a stronger understanding of those applications, but are not limited to that and can be used in several results in infinitary combinatorics (see for example \cite{KomjathBook}). One specific combinatorial property we are interested in this paper is the ultrafilter number, which has been extensively studied in recent years, as we will see in the next subsection: 
\subsection{The ultrafilter number}
The ultrafilter number for a cardinal number $\kappa$, determines how many sets one needs in order to generate an ultrafilter over $\kappa$. Let us be more precise here:
\begin{definition} Let $U$ be an ultrafilter over a cardinal $\kappa$, define:
    \begin{enumerate}
        \item a subset of an ultrafilter $U$, $\mathcal{B}\subseteq U$ is called a \textit{base} if $\forall A\in U$ there is $B\in \mathcal{B}$ such that $B\subseteq^* A$\footnote{ The order $A\subseteq^* B$ is defined by $A\setminus B$ is bounded in $\kappa$.}.
        \item The \textit{characteristics of} $U$ is $Ch(U):=\min(|\mathcal{B}|\mid \mathcal{B}\text{ is a base for }U)$
        \item The ultrafilter number $\mathfrak{u}_\kappa:=\min(Ch(U)\mid U\text{ is a uniform ultrafilter over }\kappa)$\footnote{An ultrafilter $U$ over $\kappa$ is uniform if for every $X\in U$, $|X|=\kappa$.}
    \end{enumerate}
\end{definition}
The number $\mathfrak{u}_\kappa$ is a
\textit{generalized characteristic cardinal of the continuum} as it is known that for every $\kappa$, $\kappa^+\leq \mathfrak{u}_\kappa\leq 2^\kappa$. 
As with other characteristic cardinals, the basic question is whether they can be (namely, is it consistent to) separated from the continuum. Kunen (Exercise (A10) of Chapter XIII in \cite{Kunen}) proved that using a suitable iteration (of Mathias forcing) over a model of CH, one can force a model with $2^{\aleph_0}>\mathfrak{u}_{\aleph_0}$. Kunen's method does not generalize to greater cardinals, and he asked whether it is consistent to have $\mathfrak{u}_{\aleph_1}<2^{\aleph_1}$. This question is open and has been so since the 70s.

Assuming stronger assumptions on the cardinal $\kappa$, which are known as \textit{large cardinal assumptions}, Gitik and Shelah  \cite[Lemma 1.9]{GSOnDO} forced the existence of a  cardinal $\aleph_0<\kappa$ with $2^\kappa>\mathfrak{u}_\kappa$. This cardinal $\kappa$ is extremely greater than $\omega_1$ and do not lay near the area where other mathematics occurs. The large cardinal assumption was then improved by Garti and Shelah \cite{GartiShelah2014} to a supercompact cardinal, where they used a similar iteration to the one Kunen used, but still considered an extremely large cardinal. This was further investigated by Brooke-Taylor, Fischer, Friedman, and Montoya \cite{BROOKETAYLOR201737}. Recently, a remarkable result of Raghavan and Shelah \cite{RagShel} established the consistency of $\mathfrak{u}_\kappa<2^\kappa$ for $\kappa=2^{\aleph_0}$ where they started with a much smaller large cardinal- a \textit{measurable cardinal}. However, in their model $2^{\aleph_0}$ is still very large. They also obtained the result on the much smaller cardinal $\aleph_{\omega+1}$ but starting again from a supercompact cardinal.

While all the work mentioned above concerns the ultrafilter number for regular cardinals, a different line of research about the ultrafilter number on singular cardinals has also been studied in the last decade.  
The first result in this direction uses PCF theory and is due to Shelah and Garti   \cite[Theorem 1.4]{GartiShelahUlt}:
\begin{theorem}
Suppose that $\kappa=\cf(\lambda)<\lambda$ are two cardinals, $\lambda$ is strong limit, and $\l \lambda_i\mid i<\kappa\r$ is an unbounded and increasing sequence in $\lambda$ and $E$ is an ultrafilter such that:
\begin{enumerate}
    \item $E$ is a uniform ultrafilter over $\kappa$.
    \item Each $U_i$ is a uniform ultrafilter over $\lambda_i$ carrying a strong-base $\l A_{i,\beta}\mid \beta<\theta_i\r$.
     \footnote{Strong-base for $U$ is a sequence $\l A_\alpha\mid \alpha<\theta\r$ which is $\subseteq^*$-decreasing, each $A_\alpha\in U$, and for every $B\in U$, exists $\alpha<\theta$ such that $A_\alpha\subseteq^* B$.}\end{enumerate}

     Then there is a uniform ultrafilter $U$ over $\lambda$ such that $$Ch(U)\leq tcf(\prod_{i<\kappa}\lambda_i,<_E)\cdot tcf(\prod_{i<\kappa}\theta_i,<_E).$$

\end{theorem}

This theorem provides a way to construct a uniform ultrafilter with a small base, and one can apply it in various models where these $tcf$'s have known values. For example, for any fixed $\kappa$, Garti and Shelah \cite{GartiShelahPol} have a model where $2^\lambda>\lambda^+$, a sequence of measurables $\langle \lambda_i \mid i<\kappa \rangle$, $2^{\lambda_i}=\lambda_i^+$ and $tcf(\prod_{i<\kappa} \lambda_i,<_E)=tcf(\prod_{i<\kappa}\lambda_i^+,<_E)=\lambda^+$. The measures $U_i$ can be any normal ultrafilter over $\lambda_i$, then by normality and $2^{\lambda_i}=\lambda_i^+$ when can get $\theta_i=\lambda_i^+$. So they obtained the following:
\begin{theorem}
Suppose that there is a supercompact cardinal, then there is a model with a singular cardinal $\lambda$ such that $\mathfrak{u}_\lambda=\lambda^+$ and $2^\lambda>\lambda^+$.
\end{theorem}
Then Garti, Magidor and Shelah \cite{GSM} used the single extender-based forcing to get the following:
\begin{theorem}
Suppose that $\kappa<\lambda$ are such that $\kappa$ is strong and $\lambda>\kappa$ is a limit of measurable cardinals $\l \lambda_i\mid i<\theta\r$. Then in the generic extension by the extender based forcing with $E$ being a $(\kappa,\lambda)$-extender, for every $i<\theta$, there are $\omega$-sequences of measurables $\l \lambda_{i,n}\mid n<\omega\r$ (corresponding to the measure $U_{\lambda_i})$ in $\kappa$ such that $tcf(\prod_{n<\omega}\lambda_{i,n}/J_{bd})=\lambda_i^+$. $2^\kappa\geq\lambda$ and $GCH_{<\kappa}$.  In particular, there are  ultrafilters $U_i$  over $\kappa$ such that $Ch(U_i)=\lambda_i^+$.
\end{theorem}

From the previous model we can put collapses and obtain other values of $Ch(U)$.
 

It is natural to use Magidor-Radin extender-based forcing of Merimovich \cite{Carmi4} to push these results to uncountable cofinalities. Indeed, Cummings and Morgan \cite{CummingsMorganUlt} managed to obtain it and proved the following:

\begin{theorem}
Let $\rho<\kappa<\lambda$ where $\rho$ is regular and uncountable, $\lambda$ is the least inaccessible limit of measurable cardinals greater than $\kappa$, and there is a Mitchell increasing sequence $\l E_i\mid i < \rho\r$ such that each extender $E_i$ witnesses that $\kappa$ is $\lambda$-strong and
is such that  ${}^{\kappa}Ult(V, E_i)\subseteq Ult(V, E_i)$. Then there is a cardinal-preserving generic
extension in which $cf(\kappa) = \rho$, $2^\kappa=\lambda$, and $Sp_{\chi}(\kappa)$ is unbounded in $\lambda$.
\end{theorem}
Finally, Gitik, Shelah, and Garti \cite{GGSUlt} brought everything down to $\aleph_{\omega}$ and got the result for $\aleph_{\omega}$. In the first part of this paper  we will use the forcing \cite{Sittinon} to tackle the question of separating $\mathfrak{u}_{\aleph_{\omega_1}}$ and $2^{\aleph_{\omega_1}}$ and prove the following:
\begin{theoremx}
Suppose $\kappa$ is a singular cardinal, $\rho<\kappa$ is regular and $\l \kappa_i\mid i<\rho\r$ is a sequence of strong cardinals with limit $\kappa$ and $\rho<\kappa_0$. Suppose that $\vec{E}=\l E_i \mid i<\rho\r$ is a Mitchell increasing sequence of extenders witnessing that $\kappa_i$ is $\kappa^{++}$-strong, then after forcing with $\mathbb{P}_{\vec{E}}$ we obtain a model where $\kappa=\aleph_{\omega_1}$, $2^\kappa>\aleph_{\omega_1+1}$ and $\mathfrak{u}_\kappa=\aleph_{\omega_1+1}$.
\end{theoremx}

The proof generalizes ideas similar to the one from \cite{GGSUlt}, but also applies for the countable case. Also, we provide  some missing details in the proof from \cite{GGSUlt}. 

\subsection{Overlapping Strong Extenders of long length}
As mentioned in the previous subsection, many results, both about the ultrafilter number and others, are known to hold at extremely large cardinals. One important task is to verify whether these results hold at more \textit{down to earth} cardinals. This problem is a typical problem when we force with so-called \textit{Prikry-type} forcings. The idea is to start with a large cardinal and to force a model where we destroy some properties of this large cardinal while other properties survive the forcing. This leads to solutions of many important open problems, such as the \textit{Singular Cardinal Hypothesis} \cite{MagAnnals}. This type of forcing is in extensive use in modern set theory. The major disadvantage of such a forcing is that even though the initial large cardinal loses its large cardinal property, it might still be located in a very high spot of the mathematical universe, namely above many other large cardinals, and therefore, irrelevant to solve problems at lower cardinals. Fortunately, a mechanism of crossing the gaps between the lower cardinals and the higher ones is also available in some situations. This is the so-called Prikry-type forcing with interleaved collapses, which both preserves the large cardinal properties and brings everything down to much lower cardinal (see for example, Chapter 4 of\cite{Gitik2010}). While Prikry-type forcings usually singularize a cardinal, some variations of the Radin forcing \cite{Radin} and of Gitik's overlapping extender-based forcing \cite{BenZhang} can keep $\kappa$ regular. In Section \ref{longextenderforcing}, we show how to incorporate collapses with these kinds of forcings and develop the \textit{Long Overlapping strong extender with collapses}, which is in the spirit of the one from \cite{Sittinon} or prior to that of Gitik's overlapping extender based forcing with collapses. The main innovation of our forcing is that it can keep $\kappa$ inaccessible, and by involving collapses, turn $\kappa$ into the first inaccessible. This fact is relevant for those problems in set theory seeking for consistency results at  the first inaccessible cardinal, and we believe that this forcing might be useful to tackle such problems. For example, we apply this forcing to obtain the following model:
\begin{theoremx}
 Let $\kappa$ be the least inaccessible cardinal such that there is a Mitchell increasing sequence $\langle E_i \mid i<\kappa \rangle$ witnessing that each $\kappa_i$ is $\kappa^{++}$-strong. Then it is consistent that $\kappa$ is the least inaccessible cardinal, GCH holds for every regular below $\kappa$, SCH fails for every singular below $\kappa$ and for $\lambda<\kappa$ singular, $\mathfrak{u}_\lambda<2^\lambda$.
\end{theoremx}

Convention: $p \geq q$ means $p$ is stronger than $q$. For functions $f$ and $g$ with $\dom(g) \subseteq \dom(f)$, define $f \oplus g$ ({\em $f$ overwritten by $g$}) as the function $h$ with $\dom(h)=\dom(f)$, $h(x)=g(x)$ for $x \in \dom(g)$ and $h(x)=f(x)$ otherwise.

\section{A small ultrafilter number at $\aleph_{\omega_1}$}
We start with a basic definition.
\begin{definition}
\label{ultrafilternumberdefinition}
Let $\kappa$ be an infinite cardinal and let $U$ be a uniform\footnote{A filter $F$ over a cardinal $\kappa$ is uniform iff for every $X\in F$, $|X|=\kappa$.} ultrafilter over $\kappa$. 
\begin{enumerate}
\item A {\em base} for $U$ is a collection $\mathcal{B} \subseteq U$ such that for every $B \in U$, there is $A \in \mathcal{B}$ such that $A \subseteq^* B$, namely, $A \setminus B$ is a bounded subset of $\kappa$.
\item $Ch(U):=\min\{|\mathcal{B}| \mid \mathcal{B} \subseteq U \text{ is a base for } U\}$
\item The ultrafilter number of $\kappa$, denoted by $\mathfrak{u}_\kappa$, is 
$$\min\{Ch(U) \mid U \text{ is a uniform ultrafilter over $\kappa$}\}.$$
\end{enumerate}
\end{definition}
By Claim 1.2 of \cite{GartiShelahUlt}, for any infinite $\kappa$, $\kappa<\mathfrak{u}_\kappa\leq 2^\kappa$ and if $2^\kappa=\kappa^+$, then $\mathfrak{u}_\kappa=2^\kappa$.
It is possible to have $Ch(U)<2^\kappa$ being singular \cite{Gitik2020}.
We say that $\langle W,\leq_W\rangle$ is a \textit{pre-order} if it is reflexive and transitive. The terms ``dense" and ``open" take their usual meanings.
Let us use some of the definitions of Garti, Gitik, and Shelah from \cite{GGSUlt}:
\begin{definition}
 Let $\l W_i\mid i<\lambda\r$ be a sequence of pre-orders and $F$ be a filter on $\lambda$. A \textit{Sullam} in $(\prod_{i<\lambda}W_i,F)$ is a sequence $\l f_\alpha\mid \alpha<\mu\r\subseteq \prod_{i<\lambda}W_i$ such that:
    \begin{enumerate}
        
        \item $\langle f_\alpha \mid \alpha<\mu \rangle$ is increasing $modF$, namely, if $\alpha<\beta<\mu$, then $$\{i<\lambda\mid f_\alpha(i)<_{W_i}f_\beta(i)\}\in F.$$
        \item for any list $\l V_i\mid i<\lambda\r$ such that $V_i\subseteq W_i$ is open dense, there is $\alpha<\mu$ such that
        $$\{i<\lambda\mid f_\alpha(i)\in V_i\}\in F.$$
    \end{enumerate}
 
        \end{definition}
        We only focus on  filters which extend the dual filter of the bounded ideal $J^\lambda_{bd}=\{X \subseteq \lambda \mid X\text{ is bounded in }\lambda\}$. The collection of positive sets $F^+$ has the usual meaning, namely, if $I_F=\{\lambda \setminus X \mid X \in F\}$ is the dual ideal, then $F^+:=\{X \subseteq \lambda \mid X \not \in I_F\}$.

    Let $F$ be a filter over $\lambda$ and $W$ be a pre-order. A function $g:W\rightarrow F^+$ is said to be \textit{order preserving} when for every $p \leq_{W_i} q$, $g(q) \subseteq g(p)$. We say that $g$ is \textit{deciding} if for every $A\subseteq \lambda$ and any $w\in W$, there is $w\leq_W u$ such that $g(u)\subseteq A$ or $g(u)\subseteq \lambda\setminus A$.
    \begin{definition}
     Given a singular cardinal $\lambda=\cf(\mu)<\mu$, a \textit{nice system} $\mathcal{S}$ consists of the following data:
    \begin{enumerate}
        \item A cofinal sequence $\l \lambda_i\mid i<\lambda\r$ in $\mu$ consisting of regular cardinals.
        \item A sequence $\l D_i\mid i<\lambda\r$ such that each $D_i$ is a uniform filter over $\lambda_i$.
        \item A sequence $\l W_i\mid i<\lambda\r$ of pre-orders.
        \item Functions $g_i:W_i\rightarrow D_i^+$
   which are order-preserving and deciding.
    \end{enumerate}
    \end{definition}
The following is a slight variation of \cite[Theorem 1.3]{GGSUlt}:

\begin{theorem}\label{sufficient}
    Suppose that $\lambda=\cf(\mu)<\mu$, $\mathcal{S}$ is a nice system and $D$ is a uniform ultrafilter over $\lambda$. Suppose that $\theta\in(\mu,2^\mu]$ is a regular cardinal such that:
    \begin{enumerate}
        \item $Ch(D)\leq \theta$.
        \item there is a Sullam $\l f_\beta\mid \beta<\theta\r$ in $(\prod_{i<\lambda}W_i,D)$.
        
    \end{enumerate}
    Then there is a uniform ultrafilter $U$ over $\mu$ such that $Ch(U)\leq \theta$.
\end{theorem}
\begin{proof}
    The definition of $U$ is as follows, for $X\subseteq \mu$:
    $$X\in U\Longleftrightarrow \exists \alpha<\theta \{i<\lambda\mid g_i(f_\alpha(i))\subseteq_{D_i} X\cap \lambda_i\}\in D,$$

where $\subseteq_{D_i}$ means inclusion $mod D_i$, namely $A\subseteq_{D_i} B$ if and only if $A\setminus B$ beongs to the dual ideal of $D_i$.

\begin{claim}
    $U$ is a uniform ultrafilter over $\mu$.
\end{claim}
\begin{proof}[Proof of claim.] First, note that since $\rng(g_i)\subseteq D_i^+$,  $g_i(f_\alpha(i))\neq \emptyset (mod D_i)$, hence if $X=\emptyset$, then $\{i<\lambda \mid g_i(f_\alpha(i))\subseteq_{D_i}X\cap \lambda_i\}=\emptyset$.
    Since $D$ is a uniform filter, by the definition of $U$, $\emptyset\notin U$.
    A similar argument shows that $\mu \in U$. If $X_1,X_2\in U$, then there are $\alpha_1,\alpha_2<\theta$ such that
    $$E_l:=\{i<\lambda \mid g_i(f_{\alpha_l}(i))\subseteq_{D_i}X_l\cap \lambda_i\}\in D, \ \ \text{for } l=1,2.$$
    Suppose without loss of generality that $\alpha_1<\alpha_2$. Then by the definition of Sullam, the set
    $$E_3:=\{i<\lambda\mid f_{\alpha_1}(i)<_{W_i}f_{\alpha_2}(i)\}\in D.$$
    Since all the $g_i$'s are order-preserving, for every $i\in E_3$, $g_i(f_{\alpha_2}(i))\subseteq g_i(f_{\alpha_1}(i))$. It follows that if $i\in E_1\cap E_2\cap E_3\in D$, we have that \begin{enumerate}
        \item $g_i(f_{\alpha_2}(i))\subseteq_{D_i} X_2\cap \lambda_i$ (since $i\in E_2$).
        \item $g_i(f_{\alpha_2}(i))\subseteq g_i(f_{\alpha_1}(i))\subseteq_{D_i} X_1\cap \lambda_i$ (Since $i\in E_3$ and $i\in E_1$, resp.).
    \end{enumerate}
    It follows that $g_i(f_{\alpha_2}(i))\subseteq_{D_i} X_1\cap X_2\cap\lambda_i$. By the definition of $U$, we conclude that $X_1\cap X_2\in U$.
    Showing that $U$ is closed upward is straightforward.
    To see it is an ultrafilter, let $X\subseteq \mu$. For every $i<\lambda$, consider the set $$V_i=\{ q\in W_i\mid (g_i(q)\subseteq X\cap \lambda_i)\vee (g_i(q)\subseteq \lambda_i\setminus X)\},$$ then  $V_i$ is dense. Since $g_i$ is order-preserving, $V_i$ is also open. By the definition of Sullam, there is $\alpha<\theta$ such that 
    $$F_0:=\{i<\lambda\mid f_\alpha(i)\in V_i\}\in D.$$
    This means that for each $i\in F_0$, $g_i(f_\alpha(i))\subseteq X\cap \lambda_i$ or $g_i(f_\alpha(i))\subseteq \lambda_i\setminus X$. Let us define a variable $c_i$ which in the first case above $c_i=0$ and in the second $c_i=1$. Since $D$ is an ultrafilter, there is a unique $c^*\in\{0,1\}$ such that
    $$F_1:=\{i \in F_0 \mid c_i=c^*\}\in D.$$ Suppose without loss of generality that $c^*=0$. We have that for $i \in F_1$, $g_i(f_\alpha(i)) \subseteq X \cap \lambda_i$.
    This implies that $X \in U$. 
    Finally to see that $U$ is uniform, if $X\in U$, by definition, $I=\{i<\lambda\mid X\cap \lambda_i\in D_i^+\}\in D$. Since each $D_i$ is uniform, for $i\in I$, $|X\cap\lambda_i|=\lambda_i$ and since $D$ is uniform $|X|=|\bigcup_{i\in I}X\cap\lambda_i|=\sup_{i \in I}\lambda_i=\mu$. This completes the claim.
\end{proof}
To finish the proof, let us construct a  base of size at most $\theta$ for the ultrafilter $U$. Let
\begin{enumerate}
    \item $\l d_\alpha\mid \alpha<\theta\r$ be a base for $D$ (Assumption $(1)$ of the theorem).
    \item $\l f_\beta\mid \beta<\theta\r$ be the Sullam (Assumption $(2)$ of the theorem).

\end{enumerate} 
Define for every $\alpha,\beta<\theta$, the set:
$$B_{\alpha,\beta}=\bigcup_{i\in d_\alpha}g_i(f_\beta(i)).$$

Clearly, each $B_{\alpha,\beta}$ is in $U$.
We now check that $\mathcal{B}=\{B_{\alpha,\beta}\mid \alpha,\beta<\theta\}\subseteq U$ is a base for $U$. Let $X\in U$, then by the definition, there is $\beta<\theta$ such that
$$H_0:=\{i<\lambda\mid g_i(f_{\beta}(i))\subseteq_{D_i} X\cap \lambda_i\}\in D.$$
This implies that for each $i\in H_0$, there is a set $B_i\in D_i$ such that $g_i(f_\beta(i))\cap B_i\subseteq X\cap \lambda_i$. 
 For every $i<\lambda$, define 
$$V_i:=\{q\in W_i\mid g_i(q)\subseteq B_i\vee g_i(q)\subseteq \lambda_i\setminus B_i\}.$$
This is open dense, hence by the definition of Sullam, there is $\beta'>\beta$ such that 
$$H_1=\{i<\lambda\mid f_{\beta'}(i)\in V_i\}\in D.$$
Note that if $i\in H_1$, then since $B_i\in D_i$, we have $g_i(f_{\beta'}(i))\subseteq B_i$.
Since $\beta<\beta'$, then by the definition of Sullam,
$$H_2:=\{i<\lambda\mid f_\beta(i)<_{W_i}f_{\beta'}(i)\}\in D.$$
Find $\alpha<\theta$ such that $d_\alpha\subseteq^* H_0\cap H_1\cap H_2$, and let $\zeta<\lambda$ be such that $d_{\alpha}\setminus \zeta\subseteq H_0\cap H_1\cap H_2$. To see that $B_{\alpha,\beta'}\setminus \lambda_\zeta\subseteq X$, note that if $\nu\in B_{\alpha,\beta'}\setminus \lambda_{\zeta}$, then by the definition of $B_{\alpha,\beta'}$, there is $i\in d_\alpha\setminus \zeta$ such that $\nu \in g_i(f_{\beta'}(i))$. Thus, $i\in H_0\cap H_1\cap H_2$ and therefore, $$g_i(f_{\beta'}(i))\subseteq g_i(f_{\beta}(i))\cap B_i\subseteq X\cap\lambda_i.$$
This shows that $\nu\in X$, and hence $B_{\alpha,\beta'}\subseteq^* X$.
\end{proof}
\subsection{A model where $\mathfrak{u}_{\aleph_{\omega_1}}<2^{\aleph_{\omega_1}}$}
Now let us turn to force the assumptions of Theorem \ref{sufficient} for $\mu=\aleph_{\omega_1}$ with $\theta=\aleph_{\omega_1+1}$ with $2^{\aleph_{\omega_1}}>\aleph_{\omega_1+1}$. 
\begin{notation}
\label{forcingnotation}
Fix a sequence $\l \kappa_i\mid i<\omega_1\r$. For every $\beta\leq\omega_1$ denote by $\bar{\kappa}_\beta=\sup_{\alpha<\beta}\kappa_\alpha$ and $\bar{\kappa}_0=\omega$. In particular if $\beta$ is successor then $\bar{\kappa}_{\beta}=\kappa_{\beta-1}$ and if $\beta$ is limit then $\bar{\kappa}_{\beta}<\kappa_\beta$.
 \end{notation}

Our forcing will be the one from \cite{Sittinon}, which requires the following assumptions in the ground model $V$:
\begin{itemize}
    \item GCH.
    \item a sequence $\l \kappa_i\mid i<\omega_1\r$ of strong cardinals with limit $\kappa$. In particular, $\kappa=\bar{\kappa}_{\omega_1}$ 
    \item For each $\kappa_i$, $E_i$ is a $(\kappa_i,\kappa^{++})$-extender such that $j_{E_i}:V\rightarrow M_{E_i}$ is the extender ultrapower, $M_{E_i}$ computes cardinals correctly up to and including $\kappa^{++}$, $M_{E_i}^{\kappa_i}\subseteq M_{E_i}$.
    \item For each $i$, we have $s_i:\kappa_i\rightarrow \kappa_i$ the function representing $\kappa$ in $j_{E_i}$, namely $j_{E_i}(s_i)(\kappa_i)=\kappa$. We can assume that $s_i(\nu)>\max\{\nu,\bar{\kappa}_i\}$ for every $\nu$ (see Notation \ref{forcingnotation}).
    \item For each $i_1<i_2<\omega_1$ there is a function $t_{i_2}^{i_1}: \kappa_{i_2} \to V_{\kappa_{i_2}}$ such that $j_{E_{i_2}}(t_{i_2}^{i_1})(\kappa_{i_2})=E_{i_1}$ so that $E_{i_1} \in M_{E_{i_2}}$.
    \item $\square_\kappa$ holds
\end{itemize}
The last requirement about the square will help us build a Sullam in Theorem \ref{nicesystemomega1}. The assumption can be made possible by working in some canonical model for a Woodin cardinal.

For the convenience of the reader, we include here Merimovich's notations:
\begin{itemize}
    \item For $i<\omega_1$, an \underline{$i$-domain} is a set $d\in [\kappa^{++}]^{\kappa_i}$ such that $\kappa_i+1\subseteq d$ (a set which can be the domain of the Cohen part of a condition in the extender-based forcing).
    \item Define \underline{$mc_i(d)$}=$(j_{E_i}\restriction d)^{-1}=\{\l j_{E_i}(x),x\r\mid x\in d\}$. (This is the generator of a measure used by Merimovich in his version of Extender-based forcings).
    \item Denote the measure generated by $mc_i(d)$, by \underline{$E_i(d)$}, namely $X\in E_i(d)\Longleftrightarrow mc_i(d)\in j_{E_i}(X)$.
\end{itemize}
A typical element in a measure one set of $E_i(d)$ is a sequence which provides a ``layer" of points for the continuation of the Prikry sequences appearing in the domain of a given condition. The following definition summarizes the properties we need from such sequences:
\begin{definition}
\label{object}
 An \underline{$(i,d)$-object} is a function $\mu$ such that:
 \begin{enumerate}
     \item $\kappa_i \in \dom(\mu)\subseteq d$ and  $\rng(\mu)\subseteq s_i(\mu(\kappa_i))^{++} \subseteq \kappa_i$.
      
      (Since $\dom(mc_i(d))=j_{E_i}''d$, then $j_{E_i}(\kappa_i)\in \dom(mc_i(d))\subseteq j_{E_i}(d)$. Also $\rng(mc_i(d))=d\subseteq \kappa^{++}=j_{E_i}(s_i)(\kappa_i)^{++}$).
     \item $|\dom(\mu)|=\mu(\kappa_i)<\kappa_i$ and $\mu(\kappa_i)$ is inaccessible.
     
     (since $|\dom(mc_i(d))|=|d|=\kappa_i<j_{E_i}(\kappa_i)$).
     \item $\dom(\mu)\cap \kappa_i=\mu(\kappa_i)$ and $\mu\restriction \mu(\kappa_i)=id$.
    
     (Since $\kappa_i\subseteq d$, then $\dom(mc_i(d))\cap j_{E_i}(\kappa_i)=j_{E_i}''d\cap j_{E_i}(\kappa_i)=\kappa_i$. For the second part, note that for $\alpha<\kappa_i$, $j_{E_i}(\alpha)=\alpha$ and therefore $mc_i(d)(\alpha)=\alpha$.)
     \item $\mu$ is order preserving. 
     
     (Since $j_{E_i}$ is order-preserving.)
 \end{enumerate}
 
The set \underline{$OB_i(d)$}  is the set of $(i,d)$-objects, and clearly $OB_i(d)\in E_i(d)$.
\end{definition}
We can omit the `$i$' from the ``$(i,d)$-object" and form $OB_i(d)$ since $i$ is uniquely determined by $d$ (recall that $|d|=\kappa_i$). 
\begin{definition}
If $d\subseteq d'$ are $i$-domains let $\pi_{d',d}:OB(d')\rightarrow OB(d)$ be the restriction function $\pi_{d',d}(\mu)=\mu\restriction d$ (which is equal to $\mu\restriction \dom(\mu)\cap d)$.
\end{definition}
Clearly, the generators and the measures are projected using the restriction map, namely $j_{E_i}(\pi_{d',d})(mc_i(d'))=mc_i(d)$ and $(\pi_{d',d})_*(E_i(d'))=E_i(d)$ where $(\pi_{d^\prime,d})_*$ is the natural map induced by $\pi_{d^\prime,d}$.

Here are two relevant combinatorial lemmas regarding such measures:
\begin{proposition}
 Let $0\leq i_0<i_1<...<i_n<\omega_1$ and $F:\prod_{k=0}^nA_{i_k}\rightarrow X$ is any function such that $d_{i_k}$ is $i_k$-domain, $A_{i_k}\in E_{i_k}(d_{i_k})$ and $|X|<\kappa_{i_0}$. Then there are $B_{i_k}\subseteq A_{i_k}$ such that $B_{i_k}\in E_{i_k}(d_{i_k})$ and $F\restriction \prod_{k=0}^n B_{i_k}$ is constant.
\end{proposition}
\begin{proposition} \label{boundnumberofobjects}(The bound for the number of objects with the same projection to the normal measure)
 For each $i<\omega_1$ and an 
 $i$-domain $d$, there is a set $A_i(d)$ such that
$A_i(d) \in E_i(d)$, and for each $\nu<\kappa_i$, the size of $\{\mu \in A_i(d) \mid \mu(\kappa_i) = \nu\}$ is at
most $s_i(\nu)^{++}$.
\end{proposition}
We keep the notation of $A_i(d)$. Finally, we denote the normal measure derived from $j_{E_i}$, which is $\leq_{RK}$-below $E_i$, by \underline{$E_i(\kappa_i)$}. Namely $E_i(\kappa)$ consists of all the sets $X\subseteq\kappa_i$ such that $\kappa_i\in j_{E_i}(X)$. If $A\in E_i(d)$, the projection to the normal measure is denoted by \underline{$A(\kappa_i)$} and is defined as $A(\kappa_i)=\{\mu(\kappa_i)\mid \mu\in A\}\in E_i(\kappa_i)$.
\begin{definition}
\label{forcingatomega1}
A condition in $\mathbb{P}_{\vec{E}}$ is a sequence $p=\l p_i\mid i<\omega_1\r$ such that there is a finite set $\supp(p)\in[\omega_1]^{<\omega}$, and we have that:
$$p_i=\begin{cases} \l f_i, h^0_i,h^1_i,h^2_i\r & i\in \supp(p)\\
\l f_i,A_i,H^0_i,H^1_i,H^2_i\r & i\notin\supp(p)\end{cases}$$
Such that for every every $i_1<i_2<\omega_1$, $\dom(f_{i_1})\subseteq \dom(f_{i_2})$. Denote  $\supp(p)=\{i_1<i_2<...<i_r\}$ and $i_0=0$, then for every $i<\omega_1$:
$$\bar{\kappa}_{i}<\bar{\kappa}_{i}^{+2}<f_i(\kappa_i)<s_i(f_i(\kappa_i))<s_i(f_i(\kappa_i))^+<s_i(f_i(\kappa_i))^{++}<\kappa_{i},$$
$f_i(\kappa_i)$ is inaccessible, and we require that:
\begin{enumerate}
    \item If there is $k<r$ such that $i\in [i_k,i_{k+1})$, then  $f_i$ is a partial function from $s_{i_{k+1}}(f_{i_{k+1}}(\kappa_{i_{k+1}}))^{++}$ to $\kappa_i$ such that $\kappa_i+1\subseteq\dom(f_i)$ and $|f_i|=\kappa_i$.
    \item If $i\in [i_r,\omega_1)$, then $f_i$ is a partial function from $\kappa^{++}$ to $\kappa_i$ such that $\dom(f_i)$ is an $i$-domain. In this case, abusively write the forcing in which $f_i$ lives as $Add(\kappa_i^+,\kappa^{++})$.
    \item for $i\in \supp(p)$,   $h^0_i\in Col(\bar{\kappa}_{i}^+,<f_i(\kappa_i))$, $h^1_i\in Col(f_i(\kappa_i), s_i(f_i(\kappa_i))^+)$, $h^2_i\in Col(s_i(f_i(\kappa_i))^{+3},<\kappa_i)$. 
 \item  For $i\notin\supp(p)$:
    \begin{enumerate}
        \item $A_i\in E_i(\dom(f_i))$.
        \item $\dom(H^0_i)=\dom(H^1_i)=A_i$ and $\dom(H^2_i)=A_i(\kappa_i)$.
        \item $H^0_i(\mu)\in Col(\bar{\kappa}_{i}^+,<\mu(\kappa_i))$, $H^1_i(\mu)\in Col(\mu(\kappa_i),s_i(\mu(\kappa_i))^+)$ and $H^2_i(\mu(\kappa_i))\in Col(s_i(\mu(\kappa_i))^{+3},<\kappa_i)$.
    \end{enumerate}
\end{enumerate}
\end{definition}
If $p$ is a condition, we usually represent each component of $p$ by putting the superscript $p$ to that component. For example, $f_i$ in $p$ is denoted by $f_i^p$. We also write $\dom(f_i^p)$ as $d_i^p$.
\begin{definition}
The direct order is defined by $p\leq^* q$ if $\supp(p)=\supp(q)$, for every $i$, $f^p_i\subseteq f^q_i$ and
\begin{enumerate}
\item if $i\in \supp(p)$, then for $h^{r,p}_i\leq h^{r,q}_i$ for $r=0,1,2$.
    \item if $i\notin \supp(p)$, $\pi_{\dom(f^q_i),\dom(f^p_i)}[A^q_i]\subseteq A^p_i$. $H^{r,p}_i(\pi_{\dom(f_i^q),\dom(f_i^p)}(\mu))\leq H^{r,q}_i(\mu)$ for every $\mu$ and $r=0,1$, and $H_i^{2,p}(\gamma) \leq H_i^{r,q}(\gamma)$ for every $\gamma$.
\end{enumerate}
\end{definition}
\begin{definition} Let $i\notin \supp(p)$.
$\mu\in A^p_{i}$ is {\em addable} to $p$ if:
\begin{enumerate}
    \item $\bar{\kappa}_{i}<\mu(\kappa_i)$ is inaccessible.
    \item $\cup_{\alpha<i}\dom(f_\alpha)\subseteq \dom(\mu)$ and $\mu\restriction \bar{\kappa}_{i}=id$.
    \item For every $\beta\in (\max(\supp(p)\cap i),i)$, $\{\nu\circ\mu^{-1}\mid \nu\in A^p_\beta\}\in t_i^\beta(\mu(\kappa_i))(\mu[\dom(f_\beta)])$ (recall that $t_i^\beta$ is a function such that $j_{E_i}(t^\beta_i)(\kappa_i)=E_\beta$).
\end{enumerate}
\end{definition}
\begin{remark}
    The collection of $\mu$'s which are addable to $p$ is of measure-one.
\end{remark}
\begin{definition}
\label{onestepextension}
Let $i\notin \supp(p)$,  $i_*=\max(\supp(p)\cap i)$ where $\max(\emptyset)=-1$, and $\mu\in A^p_{i}$, define $p+\mu$ as the condition $q$ such that 
$\supp(q)=\supp(p)\cup\{i\}$, and 
\begin{enumerate}
    \item For $r\in [0,i_*)\cup (i,\omega_1)$, $p_r=q_r$.
    \item For $r=i$, $f^q_i=f^p_i \oplus\mu$, $h^{0,q}_i=H^{0,p}_i(\mu)$, $h^{1,q}_i=H^{1,p}_i(\mu)$ and $h^{2,q}_i=H^{2,p}_i(\mu(\kappa_i))$.
    \item For $r\in [i_*,i)$, $j \geq 0$, $f^q_r=f^p_r\circ\mu^{-1}$ and if $r>i_*$, then $A^q_r=A^p_r\circ \mu^{-1}$, $H^{j,q}_r(\nu)=H^{j,p}_r(\nu\circ\mu)$ for $j=0,1$, and $H^{2,q}_j=H^{2,p}_j$. 
\end{enumerate}
 \end{definition}
 Define $p+\langle \mu_1 ,\cdots, \mu_n \rangle$ recursively by $p+(\langle \mu_1, \cdots, \mu_{n-1} \rangle)+\mu_n$. Define an ordering in $\mathbb{P}_{\vec{E}}$ by $p \leq q$ if $p+\vec{\mu} \leq^* q$ for some $\vec{\mu}$ ($\vec{\mu}$ could be empty).
 Sometimes, we apply an object, which does not necessarily appear in a condition, to the part of the condition below the level of the object. For example, we have a part $p \in \mathbb{P}_{\vec{E}}$, $d \supseteq d_i^p$, and $\mu \in OB_i(d)$, then $p \restriction i$ is considered as an element in $\mathbb{P}_{\vec{E} \restriction i}$, and if $t \in \mathbb{P}_{\vec{E} \restriction i}$, we denote $t_\mu$ a tuple obtained by ``squishing $t$ by $\mu$", namely we operate as in Definition \ref{onestepextension} $(1)$ for $r<i_*$ and $(3)$. Note that $t_\mu \in \mathbb{P}_{\langle t_i^\beta(\mu(\kappa_i)) \mid \beta<i \rangle}$.
 \begin{proposition}[Properties]
 \begin{enumerate}
     \item $\mathbb{P}_{\vec{E}}$ is $\kappa^{++}$-c.c.
     \item For every $p$, and for every $i \in \supp(p)$ the forcing above $p$ can be factored to a product   $$\mathbb{P}_{<i}\times Col(s_i(f_i(\kappa_i))^{+3},<\kappa_i)\times \mathbb{P}_{>i}$$
     Where $(\mathbb{P}_{>i},\leq^*)$ is a $\kappa_i^+$-closed forcing and   $|\mathbb{P}_{<i}|\leq s_i(f_i(\kappa_i))^{+2}<\kappa_i$. 
     \item $\mathbb{P}_{\vec{E}}$ has the Prikry property and the strong Prikry property (the strong Prikry property says that for every $p$ and a dense open set $D$, there is $p^* \geq^* p$ and $a \in[\omega_1]^{<\omega}$ such that for every $\vec{\mu} \in \prod_{i \in a} A_i^{p^*}$, $p^*+\vec{\mu} \in D$).
     \item Cardinals structure: In the extension, the $\kappa_i$'s are preserved and between $\kappa_i$ and $\kappa_{i+1}$ we preserve only $$\kappa_i^+<f_{i+1}(\kappa_{i+1})<s_{i+1}(f_{i+1}(\kappa_{i+1}))^{++}<s_{i+1}(f_{i+1}(\kappa_{i+1}))^{+++}$$
     In particular, for every $i\leq\omega_1$,  $\bar{\kappa}_i$ is preserved. $\kappa^+$ is preserved by the strong Prikry property, and above $\kappa^{++}$ we use the chain condition.
     \item If $\alpha<\omega_1$ is limit, in the extension, $\bar{\kappa}_\alpha=\aleph_\alpha$, $2^{\aleph_\alpha}=\aleph_{\alpha+3}$, $\kappa$ becomes $\aleph_{\omega_1}$ and $2^{\aleph_{\omega_1}}=\aleph_{\omega_1+2}$. (The mismatch for the cardinals of the powersets of singular cardinals is not a typo. In Section \ref{longextenderforcing} we will elaborate a slight modification so that in the extension, the cardinal behavior on singular cardinals will align uniformly).
     \item $\square_{\aleph_{\omega_1}}$ holds. (This is simply because we assume $\square_\kappa$ in the ground model, $\kappa$ and $\kappa^+$ are preserved in the extension, and $\kappa$ becomes $\aleph_{\omega_1}$).
       
 \end{enumerate}
  
 \end{proposition}

\begin{theorem}
\label{nicesystemomega1}
    After forcing with $\mathbb{P}_{\vec{E}}$, there is a nice system satisfying the assumption of Theorem \ref{sufficient} with $\theta=\aleph_{\omega_1+1}$.
\end{theorem}
\begin{proof}
    The proof is divided into two stages. The first stage is to build a nice system. The second stage is to find a uniform ultrafilter of small base and a Sullam.
    
    \textbf{Stage 1}: We fix any uniform ultrafilter $D$ over $\omega_1$ in the extension.
    Let us use the sequence $\lambda_i=\kappa_{i}$ for $i<\omega_1$ which consists of regular cardinals in the extension. Note that $\kappa_i$ was measurable in $V$, we fix a normal measure $D_i^\prime$ on $\kappa_i$ in $V$.
    Since the upper forcing $\mathbb{P}_{>i}$ does not add subsets to $\kappa_i$, $\kappa_i$ remains measurable after forcing with $\mathbb{P}_{>i}$ with the measure $D_i^\prime$. Also, by the small cardinality of $\mathbb{P}_{<i}$, we can lift any ultrapower embedding using a normal ultrafilter over $\kappa_i$ from $V^{\mathbb{P}_{>i}}$ to $V^{\mathbb{P}_{>i}\times\mathbb{P}_{<i}}$. Hence $\kappa_i$ remains measurable after forcing with $\mathbb{P}_{>i}\times \mathbb{P}_{<i}$. The embedding generates a normal measure extending $D_i^\prime$, and we still call the measure in the extension $D_i^\prime$. Clearly, the measurability fo $\kappa_i$ is destroyed by forcing $Col(s_i(f_i(\kappa_i))^{+3},<\kappa_i)$. However, if $D_i'\in V^{\mathbb{P}_{>i}\times\mathbb{P}_{<i}}=:V_1$ is a normal ultrafilter over $\kappa_i$, we can follow the construction in \cite[Section~17.1]{CummingsHand}:
    Let $j_{D_i'}:V_1\rightarrow M_1$ be the usual ultrapower embedding. Then, $j_{D_i'}(Col(s_i(f_i(\kappa_i))^{+3},<\kappa_i))$ is forcing equivalent to $Col(s_i(f_i(\kappa_i))^{+3},<\kappa_i)\times \mathbb{Q}$, where $\mathbb{Q}$ adds a collapsing function for every $\alpha\in[\kappa_i,j_{D'_i}(\kappa_i))$ to have cardinality $s_i(f_i(\kappa_i)))^{+3}$. We call this forcing $Col(s_i(f_i(\kappa_i))^{+3},I_i)$ where $I_i=[\kappa_i,j_{D_i^\prime}(\kappa_i))$. Over $V_1[G]$ (which is the generic extension by $\mathbb{P}_{\vec{E}}$), let $H$ be $\mathbb{Q}$-generic over $V_1[G]$, and in the model $V_1[G][H]$ we can lift $j_{D'_i}\subseteq j^*:V_1[G]\rightarrow M_1[G*H]$. Now in $V_1[G]$, we define an extension of $D_i'$:
    $$D_i:=\{X\subseteq \kappa_i\mid 0_{\mathbb{Q}}\Vdash_{\mathbb{Q}}\kappa_i\in {j}^*(X)\}$$
    Clearly, $D_i$ is uniform.
     Moreover, in $V_1[G]$ the forcing $(D_i^+,\supseteq)$ is isomorphic to $ro(\mathbb{Q})$\footnote{Where $ro(\mathbb{Q})$ is the complete boolean algebra of regular open cuts.}. In particular, since $\mathbb{Q}$ is a dense subset of $ro(\mathbb{Q})$, there is a dense embedding $g_i:\mathbb{Q}\rightarrow D_i^+$. Let $W_i=\mathbb{Q}$, then $g_i:W_i\rightarrow D_i^+$ is order preserving. To see that it is deciding, let $A\subseteq \kappa_i$ and let $q\in \mathbb{Q}$. Then either $A\cap g_i(q)\in D_i^+$ or $(\kappa_i\setminus A)\cap g_i(q)\in D_i^+$. Suppose without loss of generality that $A\cap g_i(q)\in D_i^+$. Then by density of $\mathbb{Q}$, there is $q'\geq q$ such that $g_i(q')\subseteq A\cap g_i(q)\subseteq A$. So far we have proven that $\l \kappa_i\mid i<\omega_1\r$, $\l W_i\mid i<\omega_1\r$, $\l D_i\mid i<\omega_1\r$ and $\l g_i\mid i<\omega_1\r$ forms a nice system. Note that $\mathbb{Q}$ is the collapse forcing in the sense of $V$, and $\dot{W}_i$ is decided from any $p$ with $i \in \supp(p)$, i.e. $\dot{W}_i=Col(s_i(\dot{f}_i(\kappa_i))^{+3},I_i)$.

\textbf{Stage 2:} We now prove requirements $(1)-({2})$ of Theorem \ref{sufficient}.
     $(1)$ is easy, since $\aleph_{\omega_1}$ is singular strong limit, and so $2^{\aleph_1}<\aleph_{\omega_1}$, then $Ch(D)\leq |D|<\aleph_{\omega_1+1}$. For (2) we need the following claim:

     \begin{claim}\label{keyclaim}
     For every $p \in \mathbb{P}$ and every sequence $\langle \dot{U}_i \mid i<\omega_1 \rangle$ such that $p \Vdash \dot{U}_i \subseteq \dot{W}_i$ and $\dot{U}_i$ is open dense, there is $p \leq^* p^*$ and a function $F:V_\kappa \to V_\kappa$ in $V$ such that for every generic $G$ of $\mathbb{P}$, there is a translation $F^G_i \in V[G]$ such that $F^G_i \subseteq \dot{W}_i[G]$ is open dense and is a subset of $\dot{U}_i[G]$.
     \end{claim}
     \begin{proof}
      Fix $p$ and $\dot{U}_i$ for $i<\omega_1$. Assume for simplicitiy that $p$ is pure.
      Build a $\leq^*$-increasing sequence $\langle p^i \mid i<\omega_1 \rangle$ such that for each $i$, $p^{i+1} \restriction (i+1)=p^i \restriction (i+1)$, and at each limit $\alpha$, we take $p^\alpha$ as a $\leq^*$-least upper bound of $\langle p^\beta \mid \beta<\alpha \rangle$.
      Let $p^0=p$. It remains to describe the construction at the successor stages.
      Let $i<\omega_1$ and suppose that $p^i$ is constructed.
      Write $p^i_{i+1}=\langle f,A,H^0,H^1,H^2 \rangle$.
      Let
      \begin{center}
          $\mathbb{R}^*=\{(g,\vec{r}) \in Add(\kappa_{i+1}^+,\kappa^{++}) \times \mathbb{P}_{\vec{E} \setminus (i+2)},\leq^*) \mid \dom(g) \text{ is a subset of the domains in the Cohen part of }\vec{r}\}.$
      \end{center}
      Let $N \prec H_\rho$ where $\rho$ is a sufficiently large regular cardinal, $p^i,\dot{W}_i,\dot{U}_i, \mathbb{P} \in N$, $\kappa_{i+1}+1 \subseteq N$, and ${}^{<\kappa_{i+1}}N \subseteq N$.
      Build an $\mathbb{R}^*$-increasing sequence $\{(f_\gamma,\vec{r}_\gamma) \mid \gamma<\kappa_{i+1}\}$ above $(f,p^i \setminus (i+2))$ such that every initial segment is in $N$, and for every $\mathbb{R}^*$-dense open set $D \in N$, there is $\gamma$ such that $(f_\gamma,\vec{r}_\gamma) \in D$.
      Let $f^*=\cup_\gamma f_\gamma$ and $\vec{r}$ be the $\leq^*$-least upper bound of $\langle \vec{r}_\gamma \mid \gamma<\kappa_{i+1} \rangle$.
      Then $d^*:=\dom(f^*)=N \cap \kappa^{++}$.
      Let $A^* \in E_i(d^*)$, $A^* \subseteq A_i(d^*)$, and $A^*$ projects down to a subset of $A$.
      Then $A^* \subseteq N$.
      In $N$, fix $\gamma<\kappa_{i+1}$. for each $\mu \in A^*$ with $\mu(\kappa_{i+1})=\gamma$, consider $q^\mu=\langle (p^i \restriction (i+1))_\mu,(H^0)(\mu \restriction d_i^{p^i}), (H^1)(\mu \restriction d_i^{p^i}) \rangle$. Define $D_{\mu,x}$ to be the collection of all $(h,g,\vec{r}) \in Col(s_i(\gamma)^{+3},<\kappa_i) \times \mathbb{R}^*$ such that the following set is open dense $$\Big\{(t,h^0,h^1) \geq q^\mu \mid \underset{(1)}{\underbrace{t^\frown \langle g \oplus \mu,h^0,h^1,h \rangle {}^\frown \vec{r} \Vdash x \not \in \dot{W}_i}}\text{ or }$$
      $$\underset{(2)}{
      \underbrace{t^\frown \langle g \oplus \mu,h^0,h^1,h \rangle {}^\frown \vec{r} \Vdash x \in \dot{W}_i\text{ and the condition decides some }y \geq x, \ y \in \dot{U}_i}}\Big\}.$$
      
      Let $D^\prime_\gamma=\{(g,\vec{r}) \in \mathbb{R}^* \mid \exists h(h,g,\vec{r})$ satisfies the strong Prikry property for every $D_{\mu,x}$ with $\mu(\kappa_{i+1})=\gamma\}$.
      Notice that $D_{\mu,x}$ is an open dense set which lays in $N$.
      Since the closure of the forcing  $Col(s_i(\gamma)^{+3},<\kappa_i) \times \mathbb{R}^*$ is $s_{i+1}(\gamma)^{+3}$, the number of such $\mu$ is $s_{i+1}(\gamma)^{++}$ and the number of such $x$ is $\kappa_i^+$, we have that $D_\gamma^\prime$ is an open dense set in $N$.
      By genericity, $(f^*,\vec{r}^*) \in D_\gamma^\prime$ with a witness $h=:(H^2)^*(\gamma)$.
      Let $p^{i+1}=p^i \restriction (i+1) {}^\frown \langle f^*,A^*,(H^0)^\prime,(H^1)^\prime,(H^2)^* \rangle {}^\frown \vec{r}^*$ where for $l=0,1$, $(H^l)^\prime$ is the natural map induced from $H^l$. For the rest of the proof, we denote $d_{i+1}:=\dom(f^*)$ as above.
      
      Take $p^*$ as the $\leq^*$-least upper bound for $p^i$. For each $(\mu,x)$ with $\mu \in A_{i+1}^{p^*}$, by the property of $D_{\mu \restriction d_{i+1},x}$ and the property of $p^{i+1}$, we have that for each $(\mu \restriction d_{i+1},x)$, there is a set $a_{\mu,x}:=a_{\mu \restriction d_{i+1},x} \in [\omega_1]^{<\omega}$ witnessing the strong Prikry property for $p^*$, namely for every $\vec{\tau} \in \prod_{\beta \in  a_{\mu,x}} A_\beta^{p^*}$, $((H^2)^*(\mu(\kappa_{i+1})),f_{i+1}^{p^*},(p^*+\langle \mu,\vec{\tau}\rangle) \setminus (i+2) \rangle \in D_{\mu \restriction d_{i+1},x}$.
      For each $\vec{\tau} \in \prod_{i \in a_{\mu,x}}$ maximal, fix a maximal antichain $B_{\mu,x,\vec{\tau}} \subseteq \mathbb{P}_{\langle t_{i+1}^\beta(\mu(\kappa_{i+1}) \mid \beta<i+1 \rangle} \times Col(\kappa_i^+,<\mu(\kappa_{i+1})) \times Col(\mu(\kappa_{i+1}),s_{i+1}(\mu(\kappa_{i+1}))^+)$ such that every element in $B_{\mu,x,\vec{\tau}}$ either satisfies $(1)$ or $(2)$.
      Define $F: V_\kappa \to V_\kappa$ as follows: for each $\mu,x,\vec{\tau}$ and $r \in B_{\mu,x,\vec{\tau}}$, define $F(\mu,x,\vec{\tau},r)=y$ if $(2)$ holds with the decision $y$. Otherwise, the value is $0$. 
      For other elements in $\dom(F)$ assign them the value $0$.

We now interpret $F[G]$ when $G$ is $\mathbb{P}_{\vec{E}}$-generic containing $p^*$. 
For each $x \in W_i:=\dot{W}_i[G]$, find a condition $s \in G$ above $p^*+\langle \mu, \vec{\tau} \rangle$ with $(s \restriction (i+1),(h^0_{i+1})^s,(h^1_{i+1})^s) \geq r$ for a unique $r \in B_{\mu,x,\vec{\tau}}$. Let $y_x=F(\mu,x,\vec{\tau},r)$. Define $F_i^G=\{y \mid y \geq y_x$ for some $x \in W_i\}$.
We see that $F_i^G$ is open dense and is a subset of $\dot{U}_i[G]$.
     \end{proof}
    
From the claim we get that if $\l U_i\mid i<\omega_1\r\in V[G]$ is a list of dense open subsets of $\l W_i\mid i<\omega_1\r$, then there is a function $F:V_\kappa\rightarrow V_\kappa\in V$ such that for every $i<\omega_1$, $F_i^{G}\subseteq U_i$ and $F_i^G$ is dense open. 
Since we have $GCH$ in $V$, we can enumerate by $\l g_\alpha\mid \alpha<\kappa^+\r$ all the functions $F:V_\kappa\rightarrow V_\kappa$ and in $V[G]$, denote $U_{\alpha,i}=(g_\alpha)^G_i$ and $\vec{U}_\alpha=\l U_{\alpha,i}\mid i<\omega_1\r$. Then the sequence $\l \vec{U}_\alpha\mid \alpha<\kappa^+\r$ has the following properties:
\begin{enumerate}
    \item Each $\vec{U}_\alpha$ is a sequence of dense open subsets $U_{\alpha,i}\subseteq W_i$.
    \item For every sequence $\l V_i\mid i<\omega_1\r$ of dense open subsets of $W_i$, there is some $\alpha<\kappa^+$ such that for every $i<\omega_1$, $U_{\alpha,i}\subseteq V_i$. 
\end{enumerate}
Let us now define a Sullam $\l f_\alpha\mid \alpha<\aleph_{\omega_1+1}\r$ modulo the filter of co-bounded subsets of $\omega_1$ in the generic extension. Fix a $\square_{\aleph_{\omega_1}}$-sequence $\langle C_\alpha \mid \alpha \in \lim(\aleph_{\omega_1+1}) \rangle$ such that each $C_\alpha$ has order-type below $\aleph_{\omega_1}$. Our induction hypothesis is that for each limit $\alpha$, if $i^*$ is the least such that the closure of $W_{i^*}$ is strictly greater than $ot(\lim(C_\alpha))$, the for $i \geq i^*$, $\langle f_\beta(i) \mid \beta \in \lim(C_\alpha) \cup\{\alpha\}\rangle$ is strictly increasing.
$f_0$ is random. 
Fix $f_\alpha$, let $f_{\alpha+1}$ be such that for all $i$, $f_\alpha(i)<_{W_i} f_{\alpha+1}(i)$ and $f_{\alpha+1}(i) \in U_{\alpha,i}$.
Now, assume $\alpha$ is limit. If $ot(C_\alpha)=\omega$, then let $f_\alpha(i)=\sup_{\beta \in C_\alpha} f_\beta(i)$. A straightforward argument shows that $f_\alpha$ is a $\leq^*$-upper bound of $\langle f_\beta \mid \beta<\alpha \rangle$.
Assume that $ot(C_\alpha)>\omega$.
Let $i^*$ be the least such that the closure of $W_{i^*}$ is greater than $ot(\lim(C_\alpha))$.
We divide further into two subcases. If $\lim(C_\alpha)$ is bounded in $\alpha$, then $\beta^*=\max(\lim(C_\alpha))$ exists. This only happens if $\cf(ot(C_\alpha))=\omega$ and $C_\alpha \setminus (\beta^*+1)$ has order-type $\omega$.
For this case, let $i^{**} \geq i^*$ be such that for $i\geq i^{**}$, $\langle f_\beta(i) \rangle {}^\frown \langle f_{\gamma(i)} \mid i \in C_\alpha \setminus (\beta+1) \rangle$ is strictly increasing. Define $f_\alpha(i)$ such that for $i \in [i^*,i^{**})$, $f_\alpha(i)=f_\beta(i)$, and for $i \geq i^{**}$, $f_\alpha(i)=\sup_{\gamma \in C_\alpha \setminus (\beta+1)} f_\alpha(i)$. Then $f_\alpha$ is a $\leq^*$ upper bound of $\langle f_\gamma \mid \gamma<\alpha \rangle$ and $f_\alpha$ satisfies the induction hypothesis.
We now consider the second subcase, which is the case where $\lim(C_\alpha)$ is unbounded in $\alpha$. In this case, for $i \geq i^*$, $\langle f_\beta(i) \mid \beta \in \lim(C_\alpha) \rangle$ is increasing.
Let $f_\alpha$ be such that for $i \geq i^*$, $f_\alpha(i)=\sup_{\beta \in \lim(C_\alpha)} f_\beta(i)$.
This completes the proof of Theorem \ref{nicesystemomega1}.

\end{proof}

From Theorem \ref{sufficient} and Theorem \ref{nicesystemomega1}, we conclude that
\begin{theorem}
\label{maintheorem1}
    Assume GCH, $\langle \kappa_\alpha \mid \alpha<\omega_1 \rangle$ is an increasing sequence of cardinals such that by letting $\kappa=\sup_{\alpha<\omega_1} \kappa_\alpha$,
    \begin{enumerate}
        \item for each $\alpha$, $\kappa_\alpha$ carries a $(\kappa_\alpha,\kappa^{++})$-extender $E_\alpha$.
        \item let $j_\alpha:V \to Ult(V,E_\alpha)$, then $Ult(V,E_\alpha)$ computes cardinals correctly up to and including $\kappa^{++}$.
        \item if $\beta<\alpha$, there is $t_\alpha^\beta$ such that $j_\alpha(t_\alpha^\beta)(\kappa_\alpha)=E_\beta$.
    \end{enumerate}
    Then, there is a forcing such that in a generic extension, $\mathfrak{u}_{\aleph_{\omega_1}}=\aleph_{\omega_1+1}<2^{\aleph_{\omega_1}}=\aleph_{\omega_1+2}$.
\end{theorem}

\begin{remark}
With the same argument, in the forcing extension from Theorem \ref{maintheorem1}, for any limit ordinal $\alpha<\omega_1$ we also get $\mathfrak{u}_{\aleph_\alpha}<2^{\aleph_\alpha}$. 
\end{remark}

\section{Forcing with a long sequence of overlapping extenders with collapses}
\label{longextenderforcing}
We start with a ground model $V$ and a sequence $\l \kappa_i\mid i<\kappa\rangle\in V$. \begin{notation}
for every $\beta\leq\kappa$ denote by $\bar{\kappa}_\beta=\sup_{\alpha<\beta}\kappa_\alpha$ and $\bar{\kappa}_0=\omega$. In particular if $\beta$ is successor then $\bar{\kappa}_{\beta}=\kappa_{\beta-1}$ and if $\beta$ is limit then $\bar{\kappa}_{\beta}\leq\kappa_\beta$. It follows that the sequence $\l \bar{\kappa}_\beta\mid\beta<\kappa\r$ has the same limit points as the sequence $\l \kappa_\beta\mid \beta<\kappa\r$. Namely, that for every limit $\beta\leq\kappa$, $\bar{\bar{\kappa}}_\beta=\bar{\kappa}_\beta$. \end{notation}

Assumptions in $V$:
\begin{itemize}
    \item GCH.
    \item For each $\kappa_i$, $E_i$ is a $(\kappa_i,\kappa^{++})$-extender such that $j_{E_i}:V\rightarrow M_{E_i}$ is the extender ultrapower, $M_{E_i}$ computes cardinals correctly up to and including $\kappa^{++}$, $M_{E_i}^{\kappa_i}\subseteq M_{E_i}$.
    \item For each $i$, we have $s_i:\kappa_i\rightarrow \kappa_i$ the function representing $\kappa$ in $j_{E_i}$, namely $j_{E_i}(s_i)(\kappa_i)=\kappa$. We can assume that $s_i(\nu)>\max\{\nu,\bar{\kappa}_i\}$ for every $\nu$.
    \item For each $i_1<i_2<\kappa$, there is a function $t_{i_2}^{i_1}: \kappa_{i_2} \to V_{\kappa_{i_2}}$ such that $j_{E_{i_2}}(t_{i_2}^{i_1})(\kappa_{i_2})=E_{i_1}$, and in particular $E_{i_1} \in \Ult(V,E_{i_2})$.
    \item For each $\alpha<\kappa$ limit, $\square_{ \bar{\kappa}_\alpha}$ holds.
\end{itemize}
 We need further assumptions to ensure that $\kappa$ is  the first inaccessible cardinal in the generic extension:
 \begin{itemize}
     \item $\kappa$ is inaccessible.
    \item $\kappa=\sup_{\alpha<\kappa}\kappa_i$ i.e. $\kappa=\bar{\kappa}_\kappa$.
    \item  there is no regular $\beta<\kappa$ such that $\beta=\bar{\kappa}_\beta$. \end{itemize}
\begin{remark}
    The assumption that $\kappa$ is the least inaccessible fix point of the sequence $\bar{\kappa}_i$ implies that for every limit $\alpha<\kappa$, $\bar{\kappa}_\alpha$ is singular. Just otherwise, $\bar{\kappa}_\alpha$ is also inaccessible, and hence $\alpha =\cf(\bar{\kappa}_\alpha)= \bar{\kappa}_\alpha$ contradicting the minimality of $\kappa$. It follows that for each $\beta<\kappa$, $\beta\leq \bar{\kappa}_\beta<\kappa_\beta$.
\end{remark}

Once again let us adopt Merimovich notations:
\begin{itemize}
    \item For $i<\kappa$, an \underline{$i$-domain} is a set $d\in [\kappa^{++}]^{\kappa_i}$ such that $\kappa_i+1\subseteq d$ 
    \item Define \underline{$mc_i(d)$}=$(j_{E_i}\restriction d)^{-1}=\{\langle j_{E_i}(x),x\rangle\mid x\in d\}$ (This is the generator of a measure used by Merimovich in his version of Extender-based forcings).
    \item Denote the measure generated by $mc_i(d)$, by \underline{$E_i(d)$}, namely $X\in E_i(d)\Longleftrightarrow mc_i(d)\in j_{E_i}(X)$.
\end{itemize}
We define a typical element in a measure one set of $E_i(d)$. It is a sequence which will provide a ``layer" of points for the continuation of the Prikry sequences appearing in the domain of a given condition. The proof is simply to reflect the properties of the generator $mc_i(d)$. 
\begin{definition}
\label{object2}
 An \underline{$(i,d)$-object} is a sequence/function $\mu$ such that:
 \begin{enumerate}
     \item $\kappa_i \in \dom(\mu)\subseteq d$, $\rng(\mu)\subseteq s_i(\mu(\kappa_i))^{++}$.
     \item $|\dom(\mu)|=\mu(\kappa_i)<\kappa_i$ and $\mu(\kappa_i)$ is inaccessible.
     \item $\dom(\mu)\cap \kappa_i=\mu(\kappa_i)$ and $\mu\restriction \mu(\kappa_i)=id$.
     \item $\mu$ is order preserving.
 \end{enumerate}
 
The set \underline{$OB_i(d)$}  is the set of $(i,d)$-objects, and $OB_i(d)\in E_i(d)$ (see the arguments in Definition \ref{object}).
If $d$ is clear from the context, and $\mu$ is an $(i,d)$-object, we denote $i_\mu=i$.
If $\vec{\mu}=\langle \mu_1 ,\cdots, \mu_n \rangle$ is a sequence of objects, where $i_{\mu_1}<\cdots<i_{\mu_n}$, denote $i_{\vec{\mu}}$ the ordinal $i_{\mu_n}$.
\end{definition}
We can omit the `$i$' from the ``$(i,d)$-object" and from $OB_i(d)$ since $i$ is determined by $d$ (recall that $|d|=\kappa_i$). 

\textbf{The projections:} At the price of complicating the notations what we gain is that the projections between the measures of the extender are just restriction:
\begin{definition}
If $d\subseteq d'$ are $i$-domains, let $\pi_{d',d}:OB(d')\rightarrow OB(d)$ be the restriction function $\pi_{d',d}(\mu)=\mu\restriction d$ (which is equal to $\mu\restriction \dom(\mu)\cap d)$.
\end{definition}
Clearly the generators and the measures are projected using the restriction map:
\begin{proposition}
 \begin{enumerate}
     \item $j_{E_i}(\pi_{d',d})(mc_i(d'))=mc_i(d)$.
     \item $(\pi_{d',d})_*(E_i(d'))=E_i(d)$, where $(\pi_{d^\prime,d})_*$ is the natural induced map from $\pi_{d^\prime,d}$.
 \end{enumerate}
 \end{proposition}
 \begin{proposition} (The bound for the number of objects with the same projection to the normal measure)
 \label{boundnumberobjects2}
 For each $i<\kappa$ and an 
 $i$-domain $d$, there is a set $A_i(d)$ such that
$A_i(d) \in E_i(d)$, and for each $\nu<\kappa_i$, the size of $\{\mu \in A_i(d) \mid \mu(\kappa_i) = \nu\}$ is at
most $s_i(\nu)^{++}$.
\end{proposition}
We keep the notation of $A_i(d)$. We also need notations for the normal measure:
\begin{itemize}
    \item The normal measure \underline{$E_i(\kappa_i)$} is the set of all $X\subseteq\kappa_i$ such that $\kappa_i\in j_{E_i}(X)$.
    \item If $A\in E_i(d)$ (then recall that $\kappa_i\in d$) and  the projection to normal is denoted by \underline{$A(\kappa_i)$} and is defined as $A(\kappa_i)=\{\mu(\kappa_i)\mid \mu\in A\}\in E_i(\kappa_i)$.
\end{itemize}
\begin{definition}
\label{forcinglong}
$\mathbb{P}_{\langle E_i \mid i<\kappa \rangle}$ is a sequence $p=\langle p_i\mid i<\kappa\rangle$ such that there is a finite set $\supp(p)\in[\kappa]^{<\omega}$, and we have that:
$$p_i=\begin{cases} \langle f_i, h^0_i,h^1_i\rangle & i\in \supp(p)\\
\langle f_i,A_i,H^0_i,H^1_i\rangle & i\notin\supp(p)\end{cases}$$
Such that for every $i_1<i_2<\kappa$, $\dom(f_{i_1})\subseteq \dom(f_{i_2})$. Denote  $\supp(p)=\{i_1<i_2<...<i_r\}$, then for every $i<\kappa$:
$$\bar{\kappa}_{i}<\bar{\kappa}_{i}^{+2}<f_i(\kappa_i)<s_i(f_i(\kappa_i))<s_i(f_i(\kappa_i))^+<s_i(f_i(\kappa_i))^{++}<\kappa_{i},$$
and we require that:
\begin{enumerate}
    \item If there is $k<r$ such that $i\in [i_k,i_{k+1})$ (where $i_0=0$), then $f_i$ is a partial function from $s_{i_{k+1}}(f_{i_{k+1}}(\kappa_{i_{k+1}}))^{++}$ to $\kappa_i$ such that $\kappa_i+1\subseteq\dom(f_i)$ and $|f_i|=\kappa_i$.
    \item If $i\in [i_r,\kappa)$, then $f_i$ is a partial function from $\kappa^{++}$ to $\kappa_i$ such that $\dom(f_i)$ is an $i$-domain. We will abusively write $``f_i \in Add(\kappa^{++},\kappa_i^+)"$.
    \item for $i\in \supp(p)$,   $h^0_i\in Col(\bar{\kappa}_{i}^{+2}, s_i(f_i(\kappa_i))^+)$, $h^1_i\in Col(s_i(f_i(\kappa_i))^{+3},<\kappa_i)$. (So in the generic extension $V[G]$ we will have:
            $\bar{\kappa}_{i}<(\bar{\kappa}_{i}^+)^{V[G]}=\bar{\kappa}_i^V<(\bar{\kappa}_{i}^{++})^{V[G]}=(s_i(f_i(\kappa_i))^{++})^V<(\bar{\kappa}_i^{+3})^{V[G]}=(s_i(f_i(\kappa_i))^{+3})^V<(\bar{\kappa}_{i}^{+4})^{V[G]}=\kappa_i.)$.

 \item  For $i\notin\supp(p)$:
    \begin{enumerate}
        \item if there is $k<r$ such that $i \in [i_k,i_{k+1})$, then $A_i \in t_{i_{k+1}}^i(f_{i_{k+1}}(\kappa_{i_{k+1}}))(\dom(f_i))$.
        \item if $i>i_r$, then $A_i\in E_i(\dom(f_i))$.
        \item $\dom(H^0_i)=A_i$ and $\dom(H^1_i)=A_i(\kappa_i)$.
        \item $H^0_i(\mu)\in Col(\bar{\kappa}_{i}^+,s_i(\mu(\kappa_i))^+)$ and $H^2_i(\mu(\kappa_i))\in Col(s_i(\mu(\kappa_i))^{+3},<\kappa_i)$.
    \end{enumerate}
    
\end{enumerate}
If $p$ is a condition, we usually represent each component of $p$ by putting the superscript $p$ to that component. For example, $f_i$ in $p$ is denoted by $f_i^p$. We also write $\dom(f_i^p)$ as $d_i^p$.
\end{definition}
\begin{remark}
The collapses are different from Definition \ref{forcingatomega1}. There is a flexibility to split collapses to be as in Definition \ref{forcingatomega1} or merge some collapses as in Definition \ref{forcinglong}. The reason is, we want to demonstrate a flexibility on the cardinal arithmetic on regular cardinals. Ultimately, we will obtain a ZFC model $V_\kappa$ where GCH holds at regulars, SCH fails at singulars,  small ultrafilter numbers everywhere, and $\kappa$ is the least strongly inaccessible cardinal.
\end{remark}
The direct extension is clear:
\begin{definition}
The direct order is defined by $p\leq^* q$ if $\supp(p)=\supp(q)$, for every $i$, $f^p_i\subseteq f^q_i$ and 
\begin{enumerate}
\item If $i\in \supp(p)$ then $h^{r,p}_i\leq h^{r,q}_i$ for $r=0,1$.
    \item If $i\notin \supp(p)$, $\pi_{\dom(f^q_i),\dom(f^p_i)}[A^q_i]\subseteq A^p_i$. $H^{0,p}_i(x)\leq H^{0,q}_i(\pi_{\dom(f^q_i),\dom(f^p_i)}(x))$ for all $x$, and $H_i^{1,p}(\gamma) \leq H_i^{1,q}(\gamma)$ for all $\gamma$.
\end{enumerate}
\end{definition}
\begin{definition} 
Let $i\notin \supp(p)$.
$\mu\in A^p_{i}$ is addable to $p$ if:
\begin{enumerate}
    \item $\bar{\kappa}_{i}<\mu(\kappa_i)$ is inaccessible.
    \item $\cup_{\alpha<i}\dom(f_\alpha)\subseteq \dom(\mu)$ and $\mu\restriction \bar{\kappa}_{i}=id$.
    \item For every $\beta\in (\max(\supp(p)\cap i),i), \{\nu\circ\mu^{-1}\mid \nu\in A^p_\beta\}\in t_i^\beta(\mu(\kappa_i))(\mu[\dom(f_\beta)])$.
    \end{enumerate} 
\end{definition}
\begin{remark}
    The collection of $\mu \in A_i^p$ which is addable to $p$ is of measure-one, since $i\leq \bar{\kappa}_i<\kappa_i$.
\end{remark}
\begin{definition}
\label{onestepextension2}
Let $i\notin \supp(p)$,  $i_*=\max(\supp(p)\cap i)$, where $\max(\emptyset)=-1$ and $\mu\in A^p_{i}$, define $p+\mu$ as the condition $q$ such that 
$\supp(q)=\supp(p)\cup\{i\}$, and 
\begin{enumerate}
    \item For $r\in [0,i_*)\cup (i,\kappa)$, $p_r=q_r$.
    \item For $r=i$, $f^q_i=f^p_i+\mu$, $h^{0,q}_i=H^{0,p}_i(\mu)$, and $h^{1,q}_i=H^{1,p}_i(\mu(\kappa_i))$.
    \item For $r\in [i_*,i)$ with $r \geq 0$, $f^q_r=f^p_r\circ\mu^{-1}$, if $r>i^*$, then $A^q_r=A^p_r\circ \mu^{-1}$. $H^{0,q}_r(\nu)=H^{0,q}_r(\nu\circ\mu)$  and $H^{1,q}_r=H^{1,p}_r$. Finally, if $i_* \geq 0$, then $h_{i_*}^{0,q}=h_{i_*}^{0,p}$ and  $h_{i_*}^{1,q}=h_{i_*}^{1,p}$.
\end{enumerate}
 \end{definition}

Define $p+\vec{\mu}$ recursively by $p+\langle \mu_1, \cdots, \mu_n \rangle=(p+\langle \mu_1 ,\cdots, \mu_{n-1} \rangle)+\mu_n$.
We define $p \leq q$ if $p+\vec{\mu} \leq^* q$ for some $\vec{\mu}$ ($\vec{\mu}$ could be empty).
Sometimes we apply an object which does not necessarily occur in a condition to the part of the condition lower than the level of the object, for example, we have a part $p \in \mathbb{P}_{\vec{E}}$, $d \supseteq d_i^p$, and $\mu \in OB_i(d)$, then $p \restriction i$ is considered as an element in $\mathbb{P}_{\vec{E} \restriction i}$, and if $t \in \mathbb{P}_{\vec{E} \restriction i}$, we denote $t_\mu$ a tuple obtained by ``squishing $t$ by $\mu$", namely we operate as in Definition \ref{onestepextension2} $(1)$ for $r<i_*$ and $(3)$. Note that $t_\mu \in \mathbb{P}_{\langle t_i^\beta(\mu(\kappa_i)) \mid \beta<i \rangle}$.
 
 The following definition follows from \cite{BenZhang}.
\begin{definition}
Let $p$ be a condition, $n>0$.
A $(p,n)$-fat-tree is a tree $T$ of height $n$ such that the following hold:

\begin{enumerate}
\item $Level_k(T)$ is a collection of sequences of objects of length $k+1$.
\item $Level_0(T) \in E_i(d_i^p)$ for some $i$.
\item If $\vec{\mu}=\langle \mu_1, \cdots, \mu_k \rangle \in T$ and $k<n-1$, then there is $i>i_{\vec{\mu}}$ (the definition of $i_{\vec{\mu}}$ is as in Definition \ref{object2}) such that $Succ_T(\vec{\mu}) \in E_{i}(\dom(f_i^p))$.
\end{enumerate}
We say that $T$ is {\em fully compatible} with $p$ if for every non-maximal $\vec{\mu} \in T$, $Succ_T(\vec{\mu})=A_{i_{\vec{\mu}}}^p$.

\end{definition}

The following two lemmas are easy.
\begin{lemma}
Let $p \in \mathbb{P}_{\vec{E}}$ and $T$ be a $(p,n)$-fat-tree.
\begin{enumerate}
\item If $T$ is fully compatible with $p$, then the collection $\{p+\vec{\mu} \mid \vec{\mu} \in T$ is maximal$\}$ is predense above $p$.
\item There is $p^* \geq^* p$ and a $(p^*,n)$-fat tree $T^*$ which is a subtree of $T$ of the same height and $T^*$ is fully compatible with $p^*$.
\end{enumerate}
\end{lemma}

\begin{lemma}
\label{pigeonholefattree}
    Let $T$ be a $(p,n)$-fat tree and $F:\{\vec{\mu} \in T \mid \vec{\mu}$ is maximal$\} \to \gamma$, $\gamma<\kappa_i$ where $Level_0(T) \in E_i(d)$ for some $d$.
    Then there is a fat subtree $T^\prime \subseteq T$ of the same height such that $F \restriction \{\vec{\mu} \in T^\prime \mid \vec{\mu}$ is maximal$\}$ is constant. 
\end{lemma}

\begin{lemma}[The integration lemma \cite{Sittinon}]
\label{integrationlemma}
Let $p$ be a condition, let $i \not \in \supp(p)$, $d^* \supseteq d_i^p$, $A^* \restriction d_i^p \subseteq A_i^p$. Suppose that for each $\mu \in A^*$, let $t(\mu) \geq^* (p \restriction i)_\mu$ and $h^0(\mu) \geq (H_i^0)^p(\mu \restriction d_i^p)$. Then there is $p^* \geq^* p$ such that for each $\tau \in A_i^{p^*}$ with $\mu=\tau \restriction d^*$, $(p^* \restriction i)_\tau=t(\mu)$ and $(H_i^0)^{p^*}(\tau)=h^0(\mu)$.
\end{lemma}

 \begin{theorem}[The strong Prikry property]
\label{strongtreeprikry} Let $D$ be a dense open subset of $\mathbb{P}_{\vec{E}}$ and $p$ be a condition.
Then there is a direct extension $p^* \geq^* p$ and a $(p^*,n)$-fat-tree $T$, for some $n$, which is fully compatible with $p^*$ such that for every maximal $\vec{\mu} \in T$, $p^*+\vec{\mu} \in D$.
 \end{theorem}

 \begin{remark}
     \label{technicalinduction} The proof of the strong Prikry property requires an induction on the length of the sequence of extenders. The proof where the sequence has short length was shown in \cite{Sittinon}. The proofs for the forcings on longer sequences of extenders where the lengths are below $\kappa$ are essentially the same as the proof of Theorem \ref{claimstrongprikry1}. We shall only show the strong Prikry property for $\mathbb{P}_{\vec{E}}$ while we apply the Prikry property of the forcings where the lengths of the sequences of extenders are below $\kappa$ implicitly.
 \end{remark}

\begin{proof}[Proof of Theorem \ref{strongtreeprikry}]
Let $p$ be a condition and $D$ be a dense open set.
If there is $p^* \geq^* p$ such that $p^* \in D$, then the proof is done.
Suppose it is not the case.
For simplicity, assume $p$ is pure.
The plan is to build $\langle p^n \mid 0<n<\omega \rangle$ such that for each $n$, either every direct extension of an $n$-step extension of $p^n$ is not in $D$, or there is a $(p^n,n)$-fat-tree $T^n$ fully compatible with $p^n$ such that every $n$-step extension of $p^n$ using a maximal node in $T$ is in $D$.

\underline{\textbf{Stage 1}:} \underline{step A} We build a $\leq^*$-increasing sequence $\langle q^i \mid i<\kappa \rangle$ such that
\begin{enumerate}
\item $q^0 \geq^* p$.
\item for $i^\prime<i$, $q^i_{i^\prime}=q^{i^\prime}_{i^\prime}$.
\end{enumerate}
In the end, we can take $q^*$ such that $q^*_i=q^i_i$.
Then $q^* \geq^* p$.
$q^*$ will satisfy Claim \ref{claimstrongprikry1}.
Fix $i<\kappa$.
Assume $q^{i^\prime}$ is constructed for $i^\prime<i$.
Let $q^\prime$ be such that for all $j$, $q^\prime_j$ is the weakest $``\leq^*"$-upper bound of $\{ q^{i^\prime}_j \mid i^\prime<i\}$, namely we take the union of Cohen functions, intersect the measure-one sets, and their functions whose outputs are collapses, we take the pointwise-upper bound. Note that for $i^\prime<i$, $q^\prime_{i^\prime}=q^{i^\prime}_{i^\prime}.$
Clearly, $q^\prime$ is a $\leq^*$-upper bound of $\{q^{i^\prime} \mid i^\prime<i\}$.
Write $q^\prime_i=\langle f,A,H^0,H^1 \rangle$.
Define 
\begin{center}
$\mathbb{R}_i^*=\{(g,r) \mid g\in Add(\kappa^{++},\kappa_i^+)$, $r \in (\mathbb{P}_{\langle E_\beta \mid \beta>i\rangle},\leq^*)$, and $\dom(g)$ is a subset of the domains of Cohen parts in $r\}$
\end{center}
Let $N \prec H_\theta$ for some sufficiently large $\theta$, ${}^{<\kappa_i}N \subseteq N$, $\kappa_i,q^\prime, \mathbb{P},\mathbb{R}_i^*,D \in N$, $|N|=\kappa_i$.
Enumerate dense open subsets of $\mathbb{R}_i^*$ in $N$ as $\langle D_\alpha \mid \alpha<\kappa_i \rangle$ such that every proper initial segment is in $N$.
Build an $\mathbb{R}_i^*$-increasing sequence $\langle (f_\alpha,r_\alpha) \mid \alpha<\kappa \rangle$ above $(f,q^\prime \setminus (i+1))$ such that $(f_\alpha,r_\alpha) \in D_\alpha$ for all $\alpha$.
Let $f^*=\cup_{\alpha<\kappa_i} f_\alpha$ and $r^*$ be the minimal $\leq^*$-upper bound of $\langle r_\alpha \mid \alpha<\kappa \rangle$.
Then $(f^*,r^*)$ is $(N,\mathbb{R}_i^*)$-generic in a strong sense: for $D^\prime \in N$ open dense subset of $\mathbb{R}_i^*$, there is $(f^\prime,r^\prime) \in D^\prime$ such that $(f^*,r^*) \geq (f^\prime,r^\prime) \geq (f,q^\prime \setminus (i+1))$.
Note that $d^*:=\dom(f^*)=N \cap \kappa^{++}$.
Let $A^*\in E_i(d^*)$, $A^* \subseteq A_i(d^*)$ ($A_i(d^*)$ is as in Lemma \ref{boundnumberobjects2}), and $A^*$ project down to a subset of $A$.
Then $A^* \subseteq N$.
Fix $\gamma \in A^*(\kappa_i)$.
In $N$, let $\{(t_\alpha,\mu_\alpha,h^0_\alpha) \mid \alpha<s_i(\gamma)^{++}\}$ be an enumeration of $(t,\mu,h^0)$ such that $t \in \mathbb{P}_{\langle t_i^\beta(\gamma) \mid \beta<i\rangle}$, $\mu \in A^*$ with $\mu(\kappa_i)=\gamma$, $h^0 \in Col(\bar{\kappa}_i^+,s_i(\gamma)^+)$.
Let $D_\gamma$ be the collection $(g,r) \in \mathbb{R}_i^*$ such that for all $\alpha<s_i(\gamma)^{++}$,
\begin{itemize}
\item $\dom(\mu_\alpha) \subseteq \dom(g)$.
\item there is $h^1 \geq H^1(\gamma)$ such that if there are $g^\prime \geq g \oplus \mu_\alpha$, $h^\prime \geq H^1(\gamma)$, and $r^\prime \geq^* r$ with
\begin{center}
$t_\alpha {}^\frown \langle g^\prime,h^0_\alpha,h^\prime \rangle {}^\frown r^\prime \in D$,
\end{center}
then
\begin{center}
  $t_\alpha {}^\frown \langle g \oplus \mu_\alpha,h^0_\alpha,h^1 \rangle {}^\frown r \in D$.  
\end{center}
\end{itemize}
Since $\mathbb{R}_i^*$ and $ Col(s_i(\gamma)^{+3},<\kappa_i)$ are $s_i(\gamma)^{+3}$-closed, $D_\gamma \in N$ is open dense.
By genericity, $(f^*,r^*) \in D_\gamma$ with a witness $h^1$.
Define $(H^1)^*(\gamma)=h^1$.
Let $q^i$ be such that $q^i \restriction i=q^\prime \restriction i$, $q^i_i=\langle f^*,A^*, (H^0)^*, (H^1)^* \rangle {}^\frown r^*$, where $(H^0)^*(\tau)=(H^0)(\tau \restriction \dom(f))$.
$q^i$ has the following property: for each $\mu \in A_i^{q^i}$, if $t,g,h^0,h^1$ and $r$ are such that 
\begin{center}
$t^\frown \langle g,h^0,h^1 \rangle {}^\frown r \geq q^i+\langle \mu \rangle$, $r \geq^* (q^i+\langle \mu \rangle) \setminus (i+1) ($which is $q^i \setminus (i+1))$,
\end{center}
and
\begin{center}
$t^\frown \langle g,h^0,h^1 \rangle {}^\frown r \in D$,
\end{center}
then
\begin{center}
$t^\frown \langle f^{q^i}_i \oplus \mu,h^0,(H^1)^{q^i}(\mu(\kappa_i))\rangle {}^\frown (q^i \setminus (i+1)) \in D$.    
\end{center}
Recall that we take $q^*$ such that $q^*_i=q_i^i$ for all $i$.
\begin{claim}\label{claimstrongprikry1}
For all $i<\kappa$, $\mu \in A_i^{q^*}$, if there are $t,g,h^0,h^1$, and $r$ such that
\begin{center}
$t^\frown \langle g,h^0,h^1 \rangle {}^\frown r \geq q^*+\langle \mu \rangle$,  $r \geq^* (q^*+\langle \mu \rangle) \setminus (i+1)$,
\end{center}
and
\begin{center}
$t^\frown \langle g,h^0,h^1 \rangle {}^\frown r \in D$, 
\end{center}
then
\begin{center}
$t^\frown \langle f_i^{q^*} \oplus \mu,h^0, (H_i^1)^{q^*}(\mu(\kappa_i)) \rangle {}^\frown (q^* \setminus (i+1)) \in D$.
\end{center}
\end{claim}
\begin{proof}
The claim is just a consequence of the property of $q^i$ for all $i$.
\end{proof}
\underline{Step B} We now consider the following two cases:

Case 1: For each $i<\kappa$, the collection $B_i$ of  $\mu \in A_i^{q^*}$ such that ``there are $t,h^0$ such that $t\geq^* (q^*+\langle \mu \rangle) \restriction i$, $h^0 \geq (H_i^0)^{q^*}(\mu)$, and $t {}^\frown \langle f_i^{q^*} \oplus \mu, h^0,(H_i^1)^{q^*}(\mu(\kappa_i)) \rangle {}^\frown q^* \setminus (i+1) \in D$" is of measure-zero. 
In this case, let $p^1$ be such that $p^1_i=\langle f_i^{q^*}, B_i^*, (H_i^0)^{q^*} \restriction B_i^*,(H_i^1)^{q^*} \restriction (B_i^*(\kappa_i)) \rangle$, where $B_i^*=A_i^{q^*} \setminus B_i$.

Case 2: The negation of Case 1.
This means that there is $i<\kappa$, a measure-one set $B \subseteq A_i^{q^*}$ such that for each $\mu \in B$, there are $t=t(\mu)$ and $h^0=h^0(\mu)$ such that $t \geq^* (q^*+\langle \mu \rangle) \restriction i$ and $h^0\geq (H_i^0)^{q^*}(\mu)$ such that 
\begin{center}
$t{}^\frown \langle f_i^{q^*} \oplus \mu, h^0,(H_i^1)^{q^*}(\mu(\kappa_i))\rangle {}^\frown (q^* \setminus (i+1)) \in D$.
\end{center}
Let $i$ be the least such. Use Lemma \ref{integrationlemma} with $t(\mu)$ and $h^0(\mu)$ to obtain $p^1$ such that for all $\tau \in A_i^{p^1}$ with $\mu=\tau \restriction d_i^{q^*}$, $(p^1 \restriction i)_\tau=t(\mu)$ and $(H_i^0)^{p^1}(\tau)=h^0(\mu)$.

We now have that $p^1 \geq^* p$ and for $\tau \in A_i^{p^1}$, $(p^1 \restriction i)_\tau=t(\tau \restriction d_i^{q^*})$.
With the property of $q^*$, one can check that if there is a direct extension of a one-step extension of $p^1$ entering $D$, then every one-step extension of $p^1$ using an object in $A_i^{p^1}$ is in $D$.
If this is the case, let $p^*=p^1$, and then we are done.

\underline{\textbf{Stage $n (1<n<\omega)$}:} 

\begin{remark}
\label{inductionkey}
    Note that the proof for Stage 1 holds for any condition. Furthermore, by induction hypothesis, we will assume that for any condition $r$ (in any slight variation of the long extender-based forcing with collapses, e.g. $\mathbb{P}_{\vec{E} \setminus i}$ for some $i$), and a dense open set $D^*$, for $k<n$, there is $r^* \geq^* r$ such that if there is a direct extension of a $k$-step extension of $r^*$ entering $D^*$, then there is a $(k,r^*)$-fat tree $S^*$ fully compatible with $r^*$ such that for every $\vec{\tau} \in S^*$ maximal, $r^*+\vec{\tau} \in D^*$. The statement holds for the exact proof as in Stage 1.
    \end{remark}

\underline{step A} This will be similar to step A in Stage 1, except that the dense sets we are considering in this stage are more complicated.
Suppose $p^k$ has been constructed for $k<n$.
Assume that if there is an $n-1$-step extension of $p^{n-1}$ being in $D$, then there is $(p^{n-1},n-1)$-fat tree fully compatible with $p^{n-1}$ such that every extension of $p^{n-1}$ using a maximal node in the tree belongs to $D$.
Let $p^\prime=p^{n-1}$.

Build a $\leq^*$-increasing sequence $\langle q^i \mid i<\kappa \rangle$ such that
\begin{enumerate}
    \item $p^\prime \leq^* q^0$.
    \item for $i^\prime<i$, $q_{i^\prime}^i=q_{i^\prime}^{i^\prime}$.
\end{enumerate}
Again, we take $q^*$ such that $q_i^*=q_i^i$ for all $i$ and $q^*$ will satisfy a certain property.
Fix $i<\kappa$, and assume for $i^\prime<i$, $q_{i^\prime}$ is constructed.
Let $q^\prime$ be the least $\geq^*$-upper bound of $\langle q^{i^\prime} \mid i^\prime<i\rangle$.
Hence, for $i^\prime<i$, $q^\prime_{i^\prime}=q^{i^\prime}_{i^\prime}$.
Write $q_i^\prime=\langle f,A,H^0,H^1 \rangle$, $d=\dom(f)$.
Let $\mathbb{R}_i^*$ be as in Stage 1 and
\begin{center}
$\mathbb{R}_i=\{(g,r) \mid g \in Add(\kappa^{++},\kappa_i^+), r\in (\mathbb{P}_{\langle E_\beta \mid \beta>i \rangle},\leq),$ and  \\ $\dom(g) \text{ is a subset of the domains of Cohen parts in }r\}$.
\end{center}
Note that $\mathbb{R}_i$ and $\mathbb{R}_i^*$, as sets, are equal. The difference is the ordering. Let $N \prec H_\theta$ for some sufficiently large $\theta$, ${}^{<\kappa_i}N \subseteq N$, $\kappa_i,\mathbb{P},\mathbb{R}_i^*,\mathbb{R}_i,D,q^\prime \in N$, $|N|=\kappa_i$, $\kappa_i+1 \subseteq N$, and let $\langle (f_\alpha,r_\alpha) \mid \alpha<\kappa_i \rangle$ be an $\mathbb{R}_i^*$-increasing sequence above $(f,q^\prime \setminus (i+1))$ such that every dense set contains an element in the sequence, and every proper initial segment of the sequence is in $N$. By letting $f^*=\cup_{\alpha<\kappa_i} f_\alpha$, and $r^*$ the least $\leq^*$-upper bound of $\langle r_\alpha \mid \alpha<\kappa \rangle$, then $(f^*,r^*)$ is $(N,\mathbb{R}_i^*)$-generic in the strong sense, as described in Stage 1, Step A.
Let $d^*=\dom(f^*)=N \cap \kappa^{++}$, $A^* \in E_i(d^*)$ project down to a subset of $A$ and $A^* \subseteq A_i(d^*)$.
Then $A^* \subseteq N$.

Fix $\gamma \in A^*(\kappa_i)$.
In $N$, for each $\mu$, define 
\begin{center}
    $\bar{D}_\mu=\{(g,r,h^1) \geq_{\mathbb{R}_i \times Col(s_i(\gamma)^{+3},<\kappa_i)} (f\oplus \mu, q^\prime \setminus (i+1),H^0(\gamma)) \mid$ there are $t \geq (q^\prime \restriction i)_\mu, h^0 \geq (H^0)(\mu \restriction d)$ with $t{}^\frown \langle g,h^0,h^1 \rangle {}^\frown r \in D\}$.
\end{center}
Clearly $\bar{D}_\mu \in N$ is open dense in $\mathbb{R}_i \times Col(s_i(\gamma)^{+3},<\kappa_i)$ above $(f\oplus \mu, q^\prime \setminus (i+1),H^0(\gamma))$.
We now define $D_\gamma$ as the collection of $(g,r) \in \mathbb{R}_i^*$ such that there is $(h^1)^* \geq (H_i^1)(\gamma)$ satisfying the following requirement: for each $\mu \in A^*$ with $\mu(\kappa_i)=\gamma$, 
\begin{itemize}
\item $\dom(\mu) \subseteq \dom(g)$.
\item for all $t \geq^* (q^\prime \restriction i)_\mu$, $h^0 \geq H^0(\gamma)$, if there are $g^\prime \geq g \oplus \mu$, $h^1 \geq (h^1)^*$ $a \in [\{\xi \mid i<\xi<\kappa\}]^{n-1}$, $\vec{\tau} \in \prod_{\beta \in a} A_\beta^r$, and $r^\prime \geq^* r+\vec{\tau}$ such that
\begin{center}
    $t{}^\frown \langle g^\prime, h^0,h^1\rangle {}^\frown r^\prime \in \bar{D}_\mu$,
\end{center}
then there is a $(r,n-1)$-fat tree $T$ such that for every maximal $\vec{\tau} \in T$,
\begin{center}
$t{}^\frown \langle g\oplus \mu,h,(h^1)^*\rangle {}^\frown (r+\vec{\tau}) \in \bar{D}_\mu$.
\end{center}
\end{itemize}
By our induction hypothesis as in Remark \ref{inductionkey} (we apply the remark with $\mathbb{P}_{\vec{E} \setminus (i+1)}$), the property of $p^{n-1}$, and the fact that the number of such $t$ and $h^0$ is at most $s_i(\gamma)^{+2}$, which is below the closure of $\mathbb{R}_i^*$ and $Col(s_i(\gamma)^{+3},<\kappa_i)$, we have that $D_\gamma$ is open dense in $N$. 
Hence, $(f^*,r^*) \in D_\gamma$, we obtain a witness $(h^1)^*=:(h^1)^\gamma$.
Let $(H^1)^*(\gamma)=(h^1)^\gamma$.
Let $q^i$ be such that $q^i \restriction i=q^\prime \restriction i$, $q_i^i=\langle f^*,A^*, (H^0)^*, (H^1)^* \rangle {}^\frown r^*$, where, $(H^0)^*(\tau)=(H^0)(\tau \restriction d)$.
Recall that we take $q^*$ as the least $\leq^*$-upper bound of $\langle q^i \mid i<\kappa \rangle$. We have that $q^*$ has the following property: Fix $i<\kappa$ and $\mu \in A_i^{q^*}$. Then,
\begin{itemize}
    \item either for every $t\geq (q^* \restriction i)_\mu$, $h^0 \geq (H_i^0)^{q^*}$,  $a \in [\{\xi \mid i<\xi<\kappa\}]^{n-1}$, and $\vec{\tau} \in \prod_{\beta \in a} A_\beta^{q^*}$, we have that $t {}^\frown \langle f_i^{q^*}\oplus \mu,h^0,(H_i^1)^{q^*}(\mu(\kappa_i))\rangle{}^\frown (q^* \setminus (i+1)) \not \in D$,
    \item or there are $t \geq (q^* \restriction i)_\mu$, $h^0 \geq (H_i^0)^{q^*}$, and a fat tree $T$ of height $n-1$ (not necessarily fully compatible with $q^*$) such that for every $\vec{\tau} \in T$ maximal, $t {}^\frown \langle f_i^{q^*} \oplus \mu,h^0,(H_i^1)^{q^*}(\mu(\kappa_i)) {}^\frown \vec{\tau} \in D$.
\end{itemize}
The reason that in the latter case, $T$ might not be fully compatible with $q^*$ is that by the property of $D_{\mu \restriction d_i^{q^i}}$, there is an $(n-1,q^* \setminus (i+1))$-fat tree fully compatible with $q^*$ such that for each $\vec{\tau}$ maximal in the tree, there are witnesses $t=:t_{\vec{\tau}}$ and $h^0=:h^0_{\vec{\tau}}$. We then use Lemma \ref{pigeonholefattree} to shrink the fat tree to get the fixed $t$ and $h^0$.

\underline{Step B} 
By the Prikry property (see Remark \ref{technicalinduction}) we have a possibility to choose $s \geq^* (q^* \restriction i)_\mu$ and $h^0 \geq (H_i^0)^{q^*}$ so that we have the following two cases.

\underline{Case 1}: For all $i<\kappa$, the collection $B_i$ of $\mu \in OB_i(d_i^{q^*})$ such that  ``for every  $t \geq^* (q^* \restriction i)_\mu$,  $h^0 \geq (H_0^i)^{q^*}(\mu)$, $a \in [\{\xi \mid i<\xi<\kappa\}]^{n-1}$ and $\vec{\tau} \in \prod_{\beta \in a} A_\beta^{q^*}$ with $t {}^\frown \langle f_i^{q^*} \oplus \mu,h,(H_i^i)^{q^*}(\mu(\kappa_i))\rangle {}^\frown (q^* \setminus (i+1))+\vec{\tau} \not \in D$" is of measure-one. For this case, let $p^n$ be obtained from $q^*$ by shrinking $A_i^{q^*}$ to $B_i$.

\underline{Case 2}: There is $i<\kappa$ such that the collection of $B_i$ of $\mu$ such that ``there are $t \geq^* (q^* \restriction i)_\mu$, $h^0 \geq (H_0^i)^{q^*}(\mu)$, and a $(n-1,q^* \setminus (i+1))$-fat tree $T$ such that for each $\vec{\tau} \in T$ maximal, $t(\mu){}^\frown f_i^{q^*} \oplus \mu,h^0,(H_i^1)^{q^*}(\mu(\kappa_i)) \rangle {}^\frown (q^* \setminus (i+1)) +\vec{\tau} \in D$" is of measure-one. 
Assume $i$ is the least such.
For each $\mu$, let $t=t(\mu)$ and $h^0=h^0(\mu)$ and $T=T(\mu)$ be the witnesses for the property.
Use Lemma \ref{integrationlemma} to obtain $p^n \geq^* q^*$ so that for every $\tau \in A_i^{p^n}$, $(p^n \restriction i)_\tau = t(\tau \restriction d_i^{q^*})$, and there is a $(n-1,p^n)$-fat tree $T^*(\tau)$ which projects down to a subtree of $T(\tau \restriction d_i^{q^*})$.
Let $T$ be such that $Level_0(T) \in E_i(d_i^{p^n})$, $Level_0(T)$ projects down to a subset of $B_i$, and for $\tau \in Level_0(T)$, $T_\tau=T^*(\tau)$. Shrink all relevant measure-one sets in $p^n$ so that all relevant objects appear in $T$, and finally, shrink $T$ to be fully compatible with $p^n$.

We conclude the following property of $p^n$: if there is a direct extension of an $n$-step extension of $p^n$ entering $D$, then there is a $(p^n,n)$-fat tree $T$ which is compatible with $p^n$ such that for every maximal node $\vec{\tau} \in T$, $p^n+\vec{\tau} \in D$.

Now, let $p^*$ be a $\geq^*$-upper bound of $\langle p^n \mid n<\omega \rangle$.
If $q \geq p^*$ and $q \in D$, then $q \geq^* p^*+\vec{\tau}$ for some $\rho$. Say $|\vec{\tau}|=n$. Assume that $n>1$ (the case $n=1$ is slightly simpler).
This implies that $q \geq^* p^* +\vec{\tau} \geq p^n+\vec{\tau}^\prime$, where $\vec{\tau}^\prime$ is obtained by restricting functions in $\vec{\tau}$ properly.
This implies that there is an $(p^n,n)$-tree $T$ such that for every $\vec{\mu} \in T$ maximal, $p^n+\vec{\mu} \in D$.
Let $T^*$ be a pullback of $T$ so that $T^*$ is fully compatible with $p^*$.
Then for every $\vec{\mu} \in T^*$ maximal, $p^*+\vec{\mu} \in D$.

\end{proof}

 \begin{corollary}
     $(\mathbb{P}_{\vec{E}},\leq,\leq^*)$ has the Prikry property.
     Namely for each forcing statement $\varphi$ and a condition $p$, there is $p^* \geq^* p$ such that either $p^* \Vdash \varphi$ or $p^* \Vdash \neg \varphi$.
 \end{corollary}
\begin{proof}
Let $D=\{q \mid q \Vdash \varphi$ or $q \Vdash \neg \varphi\}$.
Let $p^* \geq^* p$ and a fat tree $T$ fully compatible with $p^*$ witnessing the strong Prikry property for $D$.
By shrinking measure-one sets and the fat tree as in Lemma \ref{pigeonholefattree}, we may assume that either for all maximal $\vec{\mu} \in T$, $p+\vec{\mu} \Vdash \varphi$, or for all maximal $\vec{\mu} \in T$, $p+\vec{\mu} \Vdash \neg \varphi$.
Let $q \geq  p^*$ be such that $q$ decides $\varphi$.
Without loss of generality, assume $q \Vdash \varphi$.
Extend if necessary, assume $q \geq p^*+\vec{\mu}$ for some maximal $\vec{\mu} \in T$.
This means that for every maximal $\vec{\tau} \in T$, $p^*+\vec{\tau} \Vdash \varphi$.
Since $\{p^*+\vec{\tau} \mid \vec{\tau} \in T$ is maximal$\}$ is predense above $p^*$, we have that $p^* \Vdash \varphi$.
\end{proof}

Following the standard arguments of factorization, Prikry property, and the strong Prikry property, we have the following cardinal arithmetic

\begin{theorem}\label{maintheorem2}
After forcing with $\mathbb{P}_{\langle E_i \mid i<\kappa \rangle}$, we have that
\begin{enumerate}
    \item $\kappa$ is the first inaccessible cardinal.
    \item GCH holds for every regular cardinal below $\kappa$. Each singular cardinal below $\kappa$ is a strong limit, and SCH fails for every singular cardinal below $\kappa$.
    \item for every singular cardinal $\lambda<\kappa$, $\mathfrak{u}_\lambda=\lambda^+<\lambda^{++}=2^\lambda$.
\end{enumerate}
\end{theorem}
\begin{proof}
    \begin{enumerate}
        \item Note that $\kappa$ is strong limit. If $\kappa$ is singular, then let $p \in \mathbb{P}_{\vec{E}}$, $\alpha<\kappa$, and $\dot{f}$ be a $\mathbb{P}$-name such that $p \Vdash \dot{f}:\alpha \to \kappa$ is cofinal. Extend $p$ if necessary, assume $\alpha+1 \in \supp(p)$. Forcing above $p$ factors to $\mathbb{P}_0 \times \mathbb{P}_1$ where $|\mathbb{P}_0|=s_{\alpha+1}(f_{\alpha+1}^p(\kappa_{\alpha+1}))^{+2}$ and $(\mathbb{P}_1,\leq^*)$ is $s_{\alpha+1}(f_{\alpha+1}^p(\kappa_{\alpha+1}))^{+3}$-closed. By the Prikry property and by the closure, we can find $q \geq^* p \restriction \mathbb{P}_1$ such that for each $\gamma<\alpha$, there is a maximal antichain $A_\gamma \subseteq \mathbb{P}_0$ above $p \restriction \mathbb{P}_0$ such that for every $r \in A_\gamma$, $r^\frown q$ decides $\dot{f}(\gamma)$.  In $V$, let $X=\{\xi \mid \exists \gamma \exists r \in A_\gamma (r{}^\frown q \Vdash \dot{f}(\gamma)=\xi)\}$. Then $|X|<\kappa$ and $(p \restriction \mathbb{P}_0)^\frown q \Vdash \rng(\dot{f}) \subseteq X$, which is a contradiction.
        \item Follow the same analysis as in \cite{Sittinon}.
        \item The argument is similar to Theorem \ref{nicesystemomega1} except that we apply the version of the strong Prikry property in this section.
    \end{enumerate}
\end{proof}

\begin{remark}
    Since the forcing is $\kappa^{++}$-c.c., all cardinals above and including $\kappa^{++}$ are preserved. One can follow the argument in \cite{BenZhang} to show that $\kappa^+$ is also preserved.
\end{remark}

\begin{corollary}
    It is consistent that $GCH$ holds at every regular and $SCH$ fails at every singular $\lambda$ while $\mathfrak{u}_\lambda=\lambda^+$. 
\end{corollary}\begin{proof}
    From the previous model $V^{\mathbb{P}_{\l E_i\mid i<\kappa\r}}$, just take $M=(V^{\mathbb{P}_{\l E_i\mid i<\kappa\r}})_\kappa$ which is a $ZFC$ model (by inaccessibility of $\kappa$) which exhibit the corollary.
\end{proof}
\begin{corollary}
    It is consistent that an inaccessible $\kappa$ satisfies $\mathfrak{u}_\kappa>\kappa^+$ while there is club $C\subseteq \kappa$ such that for every $\lambda\in C$, $\mathfrak{u}_\lambda=\lambda^+<2^\lambda$.
\end{corollary}
\begin{proof}
    From the model $V^{\mathbb{P}_{\l E_i\mid i<\kappa\r}}$ force with $Add(\kappa,\kappa^{++})$, then by the $\kappa$-closure of the forcing we did not change $\mathfrak{u}_\lambda$ for $\lambda<\kappa$ and $\mathfrak{u}_\kappa=\kappa^{++}$ in the extension.
\end{proof}
\section{Open problems}
\begin{question}
What is $\mathfrak{u}_\kappa$ in the model of theorem \ref{maintheorem2}?
\end{question}

\begin{question}
\label{radin}
Can we get the same result with the Extender-based Radin-Magidor forcing with interleaved collapses?
\end{question}

We believe that the answer to Question \ref{radin} is yes, and this will reduce the consistency strength for obtaining the small ultrafilter numbers to the optimal one. 
\subsection{Acknowledgment} The authors would like to thank the referee for their careful examination of our paper and the numerous corrections and remarks they provided. Also, we would like to thank Professor Moti Gitik for his helpful remarks and long discussions.

\bibliographystyle{amsplain}
\bibliography{ref}
\end{document}